\documentclass[11pt]{amsart}

\usepackage[all]{xy}
\xyoption{dvips}
\usepackage{enumitem}
\usepackage[utf8]{inputenc}
\usepackage{pinlabel}
\usepackage{graphicx}
\usepackage{hyperref}

\newtheorem{theorem}{Theorem}[section]
\newtheorem{corollary}[theorem]{Corollary}
\newtheorem{lemma}[theorem]{Lemma}
\newtheorem{proposition}[theorem]{Proposition}

\theoremstyle{definition}

\newtheorem{example}[theorem]{Example}
\newtheorem{remark}[theorem]{Remark}

\numberwithin{equation}{section}

\newcommand{\R}{\mathbb{R}}

\newcommand{\Z}{\mathbb{Z}}
\newcommand{\C}{\mathbb{C}}

\newcommand{\id}{\mathrm{id}}

\DeclareMathOperator{\indx}{index}
\DeclareMathOperator{\coker}{coker}
\DeclareMathOperator{\Cone}{Cone}
\DeclareMathOperator{\Crit}{Crit}

\newcommand{\re}[1]{\mathfrak{Re}\,#1}
\newcommand{\im}[1]{\mathfrak{Im}\,#1}

\newcommand{\co}{\mskip0.5mu\colon\thinspace}

\title[Lifting pseudo-holomorphic discs and applications]{Lifting pseudo-holomorphic polygons to the symplectisation of $P \times \R$ and applications}
\author[G. Dimitroglou Rizell]{Georgios Dimitroglou Rizell}
\address{Department of Pure Mathematics and Mathematical Statistics, Centre for Mathematical Sciences, University of Cambridge, Wilberforce Road, Cambridge, CB3 0WB, United Kingdom}
\email{g.dimitroglou@maths.cam.ac.uk}
\thanks{This work was partially supported by the ERC Starting Grant of Fr\'{e}d\'{e}ric Bourgeois StG-239781-ContactMath.}

\begin{document}

\begin{abstract}
Let $\R \times (P \times \R)$ be the symplectisation of the contactisation of an exact symplectic manifold $P$, and let $\R \times \Lambda$ be a cylinder over a Legendrian submanifold of the contactisation. We show that a pseudo-holomorphic polygon in $P$ having boundary on the projection of $\Lambda$ can be lifted to a pseudo-holomorphic disc in the symplectisation having boundary on $\R \times \Lambda$. It follows that Legendrian contact homology may be equivalently defined by counting either of these objects. Using this result, we give a proof of Seidel's isomorphism of the linearised Legendrian contact homology induced by an exact Lagrangian filling and the singular homology of the filling.
\end{abstract}



\maketitle

\setcounter{tocdepth}{1}
\tableofcontents

\section{Introduction}

Legendrian contact homology is a Legendrian isotopy invariant which was introduced in \cite{IntroSFT} by Eliashberg, Givental, and Hofer, and in \cite{DiffAlg} by Chekanov. It associates the so-called Chekanov-Eliashberg algebra to a Legendrian submanifold $\Lambda \subset (Y,\lambda)$ of a contact manifold with a choice of contact form. This is a non-commutative differential graded algebra (DGA for short), which is freely generated by the set of Reeb chords on $\Lambda$. Roughly speaking, the differential counts rigid pseudo-holomorphic discs in a certain symplectic manifold associated to $Y$ endowed with a compatible almost complex structure, where the discs are determined by $\Lambda$ and its Reeb chords. See Section \ref{sec:background} for definitions of the above geometric objects and Section \ref{sec:background-lch} for an introduction to Legendrian contact homology.

The Chekanov-Eliashberg algebra depends on the choice of representative of the Legendrian isotopy class as well as on the compatible almost complex structure, but its homotopy type has been shown to be invariant under these choices for a wide range of ambient contact manifolds $(Y,\lambda)$. It should however be pointed out that, even though there are many contact manifolds where the invariant has been rigorously defined, there still are transversality issues that remain to be solved in order for the construction to work for a general contact manifold.

The version developed by Chekanov in \cite{DiffAlg} was the first rigorous construction of Legendrian contact homology, where it was defined for Legendrian submanifolds of the standard contact three-space
\[(\C \times \R, dz-ydx).\]
Here $z$ is a coordinate on the $\R$-factor. Legendrian contact homology was generalised to the standard contact $(2n+1)$-space 
\[\left(\C^n \times \R, dz-\sum_{i=1}^n y_idx_i\right)\]
in \cite{ContHomR} by Ekholm, Etnyre and Sullivan. The same authors \cite{ContHomP} later extended this construction to a general contactisation
\[(P \times \R,dz+\theta)\]
of an exact symplectic manifold $(P,d\theta)$, were $P$ is assumed to satisfy certain technical conditions.

The canonical projection
\[ \Pi_{\mathrm{Lag}} \co P \times \R \to P\]
of a contactisation is called the \emph{Lagrangian projection}, and will be heavily used in this paper. Observe that the Reeb chords on a Legendrian submanifold $\Lambda \subset (P \times \R,dz+\theta)$ coincide with the self-intersections of $\Pi_{\mathrm{Lag}}(\Lambda)$, and that (after a suitable perturbation) we may assume the latter to be a finite set of transverse double-points. 

Fix a compatible almost complex structure $J_P$ on $(P,d\theta)$, let $a, b_1, \hdots, b_m$ be double-points of $\Pi_{\mathrm{Lag}}(\Lambda)$, and write $\mathbf{b}=b_1\cdot \hdots \cdot b_m$. We are interested in the moduli spaces
\[ \mathcal{M}_{a;\mathbf{b}}(\Pi_{\mathrm{Lag}}(\Lambda);J_P) \]
of $J_P$-holomorphic polygons $u \co (D^2,\partial D^2) \to (P,\Pi_{\mathrm{Lag}}(\Lambda))$ having one positive boundary-puncture that is mapped to the double-point $a$, and $m$ negative boundary-punctures that are mapped to the double-points $b_1,\hdots, b_m$, where the boundary-punctures moreover appear in this cyclic order with respect to the orientation of the boundary. See Section \ref{sec:moduli-polygons} for more details.

The differential in the above versions of Legendrian contact homology is defined by counting rigid $J_P$-holomorphic polygons, that is, solutions inside a moduli space as above that moreover is required to be zero-dimensional.

The above construction of Legendrian contact homology obviously depends heavily on the existence of the projection $\Pi_{\mathrm{Lag}}$. Consequently it does not extend to a general contact manifold.

Following the philosophy of symplectic field theory \cite{IntroSFT}, relative symplectic field theory was constructed in \cite{RationalSFT} by Ekholm. This construction can be specialised to give a Legendrian contact homology for more general contact manifolds $(Y,\lambda)$, where $\lambda$ satisfies certain technical assumptions. For example, this version is well-defined and satisfy invariance for the standard contact spheres $(S^{2n+1},\lambda)$, equipped with a certain perturbation of the standard contact form. We will however only be interested in applying this theory to contactisations $(P \times \R,dz+\theta)$, where it also is well-defined.

Fix a cylindrical almost complex structure $J$ on the symplectisation $(\R \times Y,d(e^t\lambda))$ of a contact manifold $(Y,\lambda)$, where $t$ is a coordinate on the $\R$-factor. Let $a,b_1,\hdots, b_m$ be Reeb chords on $\Lambda$ and write $\mathbf{b}=b_1\cdot \hdots \cdot b_m$. We are interested in the moduli spaces
\[ \mathcal{M}_{a;\mathbf{b}}(\R \times \Lambda;J) \]
of $J$-holomorphic discs $\widetilde{u} \co (D^2,\partial D^2) \to (\R \times (P \times \R), \R \times \Lambda)$ having one positive boundary-puncture asymptotic to the Reeb chord $a$ at $t=+\infty$, and $m$ negative boundary-punctures asymptotic to the Reeb chords $b_1,\hdots,b_m$ at $t=-\infty$, where the boundary-punctures moreover appear in this cyclic order with respect to the orientation of the boundary. Observe that, since $J$ is cylindrical, these moduli-spaces have a natural $\R$-action induced by translations of the $t$-coordinate. See Section \ref{sec:moduli-discs} for more details.

The latter definition of Legendrian contact homology is equipped with a differential defined by counting non-trivial $J$-holomorphic discs (that is, not coinciding with a trivial strip $\R \times c$ over a Reeb chord) as defined above, that moreover are required to be rigid up to translation.

This version of Legendrian contact homology has the advantage that it fits more directly into the algebraic framework of symplectic field theory \cite{IntroSFT}. Indeed, it was shown in \cite{RationalSFT} that an exact Lagrangian cobordism in $\R \times Y$ between Legendrian submanifolds induce a DGA-morphism between the respective Chekanov-Eliashberg algebras.

It is a natural question whether the two versions of Legendrian contact homology defined for a contactisation, as described above, are equivalent. This has been expected and, indeed, was shown to be true for the standard contact three-space $\C \times \R$ in \cite[Theorem 7.7]{InvLegCoh}.

Also, in \cite[Lemma 2.13]{NonIsoLeg} a generalisation of this result to the general contact $(2n+1)$-space $\C^n \times \R$ was outlined. One goal of this paper is to fill in the technical details of this proof, and also to obtain the result for Legendrian contact homology in general contactisations.

\section*{Acknowledgements}
The author would like to thank Frédéric Bourgeois, Tobias Ekholm, Roman Golovko, and Samuel Lisi for valuable discussions and for showing interest in this work. The author is also grateful to Baptiste Chantraine and Paolo Ghiggini for good feedback and for pointing out the possibility of a certain degeneration of pseudo-holomorphic discs. Finally, the author would like to thank to the referee, who pointed out a gap.

\section{Results}
In the following we will always assume that $J_P$ is a compatible almost complex structure on $(P,d\theta)$ for which $(P,d\theta,J_P)$ has \emph{finite geometry at infinity}, see \cite[Definition 2.1]{ContHomP}. For instance, if $(P,d\theta)$ is symplectomorphic to the completion of a Liouville domain, such almost complex structures always exist. Also, this is true for the cotangent bundle of a smooth (not necessarily closed) manifold with its standard symplectic structure.

The main result is the following. Suppose that we are given a compatible almost complex structure $J_P$ on $(P,d\theta)$ and a Legendrian submanifold $\Lambda \subset (P \times \R,dz+\theta)$ satisfying some technical, not too restrictive, assumptions. It follows that $J_P$-holomorphic polygons in $P$ with boundary on $\Pi_{\mathrm{Lag}}(\Lambda)$ lift to pseudo-holomorphic discs in the symplectisation $\R \times (P \times \R)$ with boundary on $\R \times \Lambda$. As a direct consequence, it will follow that the two versions of Legendrian contact homology discussed above are equivalent.

An \emph{exact Lagrangian filling} (inside the symplectisation) of a Legendrian submanifold $\Lambda \subset \R \times P$ is an exact Lagrangian submanifold
\[L \subset (\R \times (P \times \R),d(e^t(dz+\theta)))\] that coincides with the cylinder $(N,+\infty) \times \Lambda$ outside of a compact set.

We use the above result for computing the wrapped Floer homology of the pair consisting of an exact Lagrangian filling and a small push-off of itself. This computation is then used to deduce properties of the linearised Legendrian contact homology of $\Lambda$. In particular, we will prove that the linearised Legendrian contact cohomology induced by the filling is isomorphic to the singular homology of the filling, a result that was first observed by Seidel.

A proof of Seidel's isomorphism was outlined in \cite{RationalSFT2}, but there it depends on a conjectural analytical lemma. We will prove this theorem following a similar strategy, but we will circumvent the analytical difficulties of the latter lemma by establishing only its algebraic consequences.

\subsection{Lifting pseudo-holomorphic polygons to the symplectisation}
\label{sec:res-lift}
A compatible almost complex structure $J_P$ on $P$ lifts to a unique cylindrical (see Section \ref{sec:cylacs}) almost complex structure $\widetilde{J}_P$ on the symplectisation
\[(\R \times (P \times \R),d(e^t(dz+\theta)))\]
determined by the requirement that the canonical projection
\[\pi_P \co \R \times (P \times \R) \to P \]
is $(\widetilde{J}_P,J_P)$-holomorphic, that is, $(D\pi_P)\widetilde{J}_P=J_P(D\pi_P)$. 

Let $\Lambda \subset P \times \R$ be a fixed chord-generic closed, not necessarily connected, Legendrian submanifold. We will be interested in almost complex structures $J_P$ that are \emph{regular} for the moduli spaces $\mathcal{M}_{a;\mathbf{b}}(\Pi_{\mathrm{Lag}}(\Lambda);J_P)$, i.e.~chosen so that the latter moduli spaces are transversely cut out.

In order to guarantee the existence of regular almost complex structures, the following real-analyticity condition will be used repeatedly.
\begin{enumerate}[label={(RA)}, ref={(RA)}]
\item \emph{There is a neighbourhood $U \subset P$ of the double-points $\mathcal{Q}(\Lambda)$ of $\Pi_{\mathrm{Lag}}(\Lambda)$ in which $J_P$ is integrable and $\Pi_{\mathrm{Lag}}(\Lambda)$ is real-analytic.} \label{ra}
\end{enumerate}
To that end, by \cite[Proposition 2.3, Lemma 4.5]{ContHomP} it follows that a compatible almost complex structure on $P$ satisfying \ref{ra} can be made regular after an arbitrarily small perturbation having compact support inside $U \setminus \mathcal{Q}(\Lambda)$.

Observe that if $J_P$ is integrable in some neighbourhood of the double-points of $\Pi_{\mathrm{Lag}}(\Lambda)$, one can perturb $\Lambda$ by a Legendrian isotopy to make its projection real-analytic in some, possibly smaller, neighbourhood containing the double-points. Alternatively, given any Legendrian submanifold $\Lambda \subset P \times \R$, it is always possible to construct a compatible almost complex structure that satisfies \ref{ra}.

\begin{theorem}
\label{thm:lift}
Let $J_P$ be a regular compatible almost complex structure on $(P,d\theta)$ that is integrable in some neighbourhood of the double-points of the generic Lagrangian immersion $\Pi_{\mathrm{Lag}}(\Lambda)$, and let $\widetilde{J}_P$ be the cylindrical lift of $J_P$ to the symplectisation of $(P \times \R,dz+\theta)$. It follows that $\pi_P$ induces a diffeomorphism
\begin{gather*}\mathcal{M}_{a;\mathbf{b}}(\R \times \Lambda;\widetilde{J}_P)/\R \to \mathcal{M}_{a;\mathbf{b}}(\Pi_{\mathrm{Lag}}(\Lambda);J_P),\\
\widetilde{u} \mapsto \pi_P \circ \widetilde{u},
\end{gather*}
where $\widetilde{J}_P$ moreover is regular for the moduli spaces $\mathcal{M}_{a;\mathbf{b}}(\R \times \Lambda;\widetilde{J}_P)$.
\end{theorem}
We refer to Section \ref{sec:proofs} for the proof.
\begin{remark}
In the case when
\[(P=\C,\theta=xdy=-(1/2)d(x^2)(J_P\cdot),J_P=i),\]
the fact that the above map is a bijection follows from \cite[Theorem 7.7]{InvLegCoh}. We generalise this result to triples $(P,\theta=-d\alpha(J_P\cdot),J_P)$ where $\alpha \co P \to \R$ is a smooth and strictly $J_P$-convex function. See Lemma \ref{lem:lift} together with Remark \ref{rem:lift} below. Note that Theorem \ref{thm:lift} is stronger, since it also asserts that the lifted cylindrical almost complex structure is regular.
\end{remark}

By the above theorem, the $J_P$-holomorphic polygons in the definition of the Legendrian contact homology differential given in \cite{ContHomP} correspond bijectively to the $\widetilde{J}_P$-holomorphic discs in the definition of the differential given in \cite{RationalSFT2}. The following corollary thus immediately follows.
\begin{corollary}
\label{cor:main}
Let $\Lambda \subset P \times \R$ be a chord-generic Legendrian submanifold of a contactisation. There are choices of regular compatible almost complex structures for which the two versions of the Chekanov-Eliashberg algebra (with coefficients in $\Z_2$) as defined in \cite{ContHomP} and \cite{IntroSFT}, respectively, become equal.
\end{corollary}

\begin{remark}
Coherent orientations for the moduli spaces of $J_P$-holomorphic polygons in exact symplectic manifolds with boundary on $\Pi_{\mathrm{Lag}}(\Lambda)$ were defined in \cite{OrientLeg}, under the additional assumption that $\Lambda$ is spin. The version of the Chekanov-Eliashberg algebra in \cite{ContHomP} can thus be defined over $\Z$ in this case. However, coherent orientations have not yet been worked out in detail for the corresponding moduli spaces in the symplectisation. This is the reason why we only acquire the result for coefficients in $\Z_2$. It is of course expected that the latter moduli spaces can be coherently oriented as well, and that the diffeomorphism in Theorem \ref{thm:lift} is orientation preserving.
\end{remark}

\subsection{Applications}
For simplicity we will only consider the case when the symplectic manifold $(P,d\theta)$ has vanishing Chern class and $\Lambda \subset P \times \R$ has vanishing Maslov class. We will furthermore assume that the Lagrangian fillings considered have vanishing Maslov class. In this case all gradings can be taken in $\Z$. The below results hold in the general situation as well, but the grading has then to be chosen in an appropriate group $\Z_\mu$.

An \emph{augmentation} of a unital DGA $(\mathcal{A}_\bullet,\partial)$ over $\Z_2$ is a unital DGA-morphism
\[\epsilon \co (\mathcal{A}_\bullet,\partial) \to (\Z_2,0).\]
Augmentations can be used to define the so called linearisation of the DGA, which is a complex over $\Z_2$ spanned by the generators. See Section \ref{sec:linearised} for more details. If the Chekanov-Eliashberg algebra has an augmentation, then we use
\[CL^\bullet(\Lambda)=(\Z_2\langle \mathcal{Q} \rangle,d_\epsilon)\]
to denote the corresponding linearised co-complex, where $\mathcal{Q}(\Lambda)$ denotes the set of Reeb chords on $\Lambda$. We will use $HCL^\bullet(\Lambda;\epsilon)$ to denote the homology of this co-complex.

Observe that the cohomology depends on the choice of augmentation, but that the set of isomorphism classes of all linearised Legendrian contact homologies is a Legendrian isotopy invariant (also, see Remark \ref{rem:stabletame}).

It is shown in \cite[Theorem 1.1]{RationalSFT2} that, in accordance with the principles of symplectic field theory, an exact Lagrangian filling of $\Lambda$ induces an augmentation of its Chekanov-Eliashberg algebra. For augmentations arising this way, there are some strong consequences for the corresponding linearised Legendrian contact cohomology. Namely, by using Theorem \ref{thm:lift} we prove Theorem \ref{thm:twocopy}, of which the following is a direct consequence.

\begin{theorem}[Seidel, Conjecture 1.2 in \cite{RationalSFT2}]
\label{thm:main}
Let $\Lambda \subset P \times \R$ be a closed Legendrian $n$-dimensional submanifold which has an exact Lagrangian filling $L \subset \R \times (P\times \R)$ inside the symplectisation. It follows that there is an isomorphism
\[\delta \co H_{n-\bullet}(L) \to HCL^{\bullet} (\Lambda;\epsilon),\]
where $\epsilon$ is the augmentation induced by the filling.
\end{theorem}

The proof, the idea of which goes back to Seidel, was given in \cite{RationalSFT} but there depends on the conjectural \cite[Lemma 4.11]{RationalSFT2}. The idea is to use the relation between the wrapped Floer homology of an exact Lagrangian filling and the linearised Legendrian contact cohomology of its end (also, see \cite{OnWrapped}), together with the observation that wrapped Floer homology vanishes inside symplectisations of contactisations (see Proposition \ref{prop:vanish}).

Wrapped Floer homology is a version of Lagrangian intersection Floer homology, originally defined in \cite{MorseTheoryLagr} by Floer,  generalised to non-compact exact Lagrangian submanifold inside a Liouville domain. It first appeared in the literature in \cite{AbbonFloer}. Different versions have later been developed in \cite{OpenStringAnalogue}, \cite{OnWrapped}, \cite{SympGeomCot} and \cite{RationalSFT2}.  We will be using the latter version, which is outlined in Section \ref{sec:background-floer}.

The below corollary of \cite[Conjecture 1.2]{RationalSFT2} was also shown in \cite{RationalSFT2}, but its proof again depends on the the conjectural \cite[Lemma 4.11]{RationalSFT2}. We will not give a proof of this lemma, but Theorem \ref{thm:twocopy} below establishes a result that gives the same information on the algebraic level. See Remark \ref{rem:conj} for a further discussion. This result is sufficient in order to establish
\begin{corollary}[Corollary 1.3 in \cite{RationalSFT2}]
\label{cor:dual}
Let $\Lambda \subset P \times \R$ be a closed Legendrian $n$-dimensional submanifold which has an exact Lagrangian filling $L \subset \R \times (P\times \R)$. The diagram below is commutative, where the horizontal sequences are exact, the upper sequence is the long exact sequence for the singular homology of a pair, and where the vertical arrows are isomorphisms.
\[
\xymatrix{
\ar[r] & H_{k+1}(\Lambda) \ar[r] \ar[d]^{\id} & H_{k+1}(L) \ar[r] \ar[d]^{\delta} & H_{k+1}(L,\Lambda) \ar[r] \ar[d]^{H^{-1} \circ \delta'} & H_k(\Lambda) \ar[d]^{\id} \ar[r] &  \\
\ar[r] & H_{k+1}(\Lambda) \ar[r]^-{\sigma} & HCL^{n-k-1}(\Lambda;\epsilon) \ar[r]& HCL_k(\Lambda;\epsilon) \ar[r]^-{\rho} & H_k(\Lambda) \ar[r] &  
}
\]
Here $H:=\begin{pmatrix}\rho & \eta\end{pmatrix}^t$, $\delta:=\begin{pmatrix} g & \sigma \end{pmatrix}$, $\delta':=\begin{pmatrix} \gamma & g \end{pmatrix}^t$ where we refer to Section \ref{sec:twocopy} for the definition of $\rho$, $\eta$, and $\sigma$, and to Theorem \ref{thm:twocopy} for the definition of $g$ and $\gamma$.
\end{corollary}
The long exact sequence at the bottom is the one constructed in \cite[Theorem 1.1]{DualityLeg}, where we have used Theorem \ref{thm:lift} above to translate these results to the Chekanov-Eliashberg algebra defined in terms of the symplectisation (see Section \ref{sec:twocopy}).
\begin{remark}
The original formulation of \cite[Theorem 1.1]{DualityLeg} requires $\Lambda \subset P \times \R$ to be \emph{horizontally displaceable}, i.e.~that $\Pi_{\mathrm{Lag}}(\Lambda) \subset P$ can be displaced from itself by a Hamiltonian isotopy. However, Proposition \ref{prop:vanishlch} shows that the results in \cite{DualityLeg} also apply for the linearised Legendrian contact homology induced by an exact Lagrangian filling in $\R \times (P \times \R)$. The reason is that such a filling always is displaceable, as follows from the proof of Proposition \ref{prop:vanish}. Observe that this holds even in the case when the Legendrian submanifold itself is not horizontally displaceable (see Example \ref{ex:nondisp}).
\end{remark}

Finally we note that results analogous to Theorem \ref{thm:main} and Corollary \ref{cor:dual} for the generating family homology have been obtained in \cite{ObsLagCob} under the assumption that the Legendrian submanifold and its filling posses generating families.

\section{General definitions}
\label{sec:background}
\subsection{Symplectic and contact manifolds}
A \emph{contact manifold} $(Y,\xi)$ is a smooth $(2n+1)$-dimensional manifold $Y$ together with a maximally non-integrable field of tangent hyperplanes $\xi$. We will consider the case when $\xi=\ker \lambda$ for a fixed one-form $\lambda$, the so called \emph{contact form}. Maximal non-integrability in this case means that the $(2n+1)$-form
\[\lambda \wedge (d\lambda)^{\wedge n} \neq 0\]
is nowhere vanishing.

An $n$-dimensional submanifold $\Lambda \subset (Y,\xi)$ of a $(2n+1)$-dimensional contact manifold is called \emph{Legendrian} if it is tangent to $\xi$. Two Legendrian submanifolds are \emph{Legendrian isotopic} if they are smoothly isotopic through Legendrian submanifolds. Determining Legendrian isotopy classes is an important, but subtle, question in contact geometry.

A choice of a contact form $\lambda$ on $Y$ determines the so called \emph{Reeb vector field} $R$ by
\[\lambda(R)=1, \:\: \iota_R d\lambda=0.\]
An integral-curve of $R$ starting and ending on two different sheets of $\Lambda$ is called a \emph{Reeb chord}.

A \emph{symplectic manifold} $(X,\omega)$ is an even-dimensional manifold together with a closed non-degenerate two-form. An $n$-dimensional submanifold $L \subset (X,\omega)$ of a $2n$-dimensional symplectic manifold is called \emph{Lagrangian} if $\omega|_{TL}=0$.

We say that the symplectic manifold is \emph{exact} if $\omega$ is exact. Observe that an exact symplectic manifold never is closed. An immersion of an $n$-dimensional manifold $L$ into an exact symplectic $2n$-manifold $(X,d\theta)$ is called \emph{exact Lagrangian} if the pull-back of $\theta$ to $L$ is an exact one-form.

The \emph{contactisation} of an exact  symplectic manifold $(P,d\theta)$ is the contact manifold
\[(P \times \R, dz+\theta),\]
where $z$ is a coordinate on the $\R$-factor. For this choice of contact form, the Reeb vector field is given by $R=\partial_z$.

The canonical projection
\[\Pi_{\mathrm{Lag}} \co P \times \R \to P\]
is called the \emph{Lagrangian projection} and it is easily checked that if $\Lambda \subset P \times \R$ is Legendrian, then $\Pi_{\mathrm{Lag}}(\Lambda) \subset (P,d\theta)$ is an exact Lagrangian immersion. Observe that the self-intersections of $\Pi_{\mathrm{Lag}}(\Lambda)$ correspond bijectively to the Reeb chords on $\Lambda$.

\begin{example}
The one-jet space of a smooth manifold $M$ can be endowed with a natural contact form
\[(J^1(M)=T^*M \times \R,dz+\theta_M),\]
where $\theta_M$ is minus the canonical one-form (also often called the Liouville form). This is also the archetypal example of a contactisation. The canonical one-form is given by $\sum_iy_idx^i$ in local canonical coordinates on $T^*M$, i.e. in coordinates of the form $(x_i,y_idx_i)$ for some choice of local coordinates $x_i$ on $M$. Specialising to the case $M=\R^n$ we obtain the standard contact $(2n+1)$-space. 
\end{example}

\subsection{Hamiltonian isotopies}
A smooth time-dependent Hamiltonian
\[H_s \co X \to \R\]
on a symplectic manifold $(X,\omega)$ gives rise to the corresponding Hamiltonian vector field $X_{H_s}$ on $X$, which is determined by
\[\iota_{X_{H_s}}\omega=dH_s.\]
We denote the induced one-parameter flow by
\begin{gather*}
\R \times X \to X,\\
(s,x) \mapsto \phi^s_H(x),
\end{gather*}
which can be seen to preserve the symplectic form. A \emph{Hamiltonian isotopy} is an isotopy of a symplectic manifold induced by a Hamiltonian as above.

It follows by Weinstein's Lagrangian neighbourhood theorem that any one-parameter family of exact Lagrangian embeddings can be realised by a time-dependent Hamiltonian isotopy of the ambient symplectic manifold.

\subsection{Exact Lagrangian cobordisms and fillings}
\label{sec:defexact}

The symplectisation of a contact manifold $(Y,\lambda)$ is the exact symplectic manifold 
\[(\R \times Y,d(e^t\lambda)),\]
where $t$ is a coordinate on the $\R$-factor. It is easily checked that $\Lambda \subset (Y,\lambda)$ is Legendrian if and only if the cylinder $\R \times \Lambda \subset \R \times Y$ is Lagrangian.

An \emph{exact Lagrangian cobordism} $L$ (inside the symplectisation) from the Legendrian submanifold $\Lambda_-$ to $\Lambda_+$ is an exact Lagrangian submanifold in $\R \times Y$ which coincides with
\[((-\infty,-N) \times \Lambda_- ) \cup ((N,+\infty) \times \Lambda_+), \:\: N > 0,\]
outside of a compact subset $L \cap \{t \in [-N,N]\}$. We moreover require there to exist a primitive of the pull-back of $e^t \lambda$ to $L$ that is constant when restricted to either of the ends $L \cap \{t \le -N\}$ or $L \cap \{t \ge N\}$ (this is automatically satisfied if the end is connected). In the case $\Lambda_-=\emptyset$, we say that $L$ is an \emph{exact Lagrangian filling} of $\Lambda_+$.

Suppose that $V$ is an exact Lagrangian cobordism from $\Lambda_-$ to $\Lambda_0$ and that $W$ is an exact Lagrangian cobordism from $\Lambda_0$ to $\Lambda_+$. We also assume that $V$ and $W$ have been translated appropriately in the $t$-direction, in order for
\begin{gather*}
V \cap \{ t \ge -1 \} = [-1,+\infty) \times \Lambda_0,\\
W \cap \{t \leq 1 \} = (-\infty,1] \times \Lambda_0,
\end{gather*}
to hold. We may then define the \emph{concatenation} of $V$ and $W$ to be
\[ V \odot W := (\{t \le 0 \} \cap V ) \cup (\{t \ge 0 \} \cap W),\]
which is a Lagrangian cobordism from $\Lambda_-$ to $\Lambda_+$. In the case when $V$ and $W$ are exact, it follows that the concatenation is exact as well. 

Consider the translation
\begin{gather*}
\tau_s \co \R \times Y \to \R \times Y,\\
\tau_s(t,y) =(s+t,y).
\end{gather*}
For each $s \geq 0$ we also define the concatenation
\[ V \odot_s W := (\{t \le 0 \} \cap V ) \cup (\{t \ge 0 \} \cap \tau_s(W)),\]
which by construction  is cylindrical in the set $\{ -1 \le t \le 1+s \}$.

Observe that all the concatenations $V \odot_s W$ are Hamiltonian isotopic.

\subsection{Almost complex structures}
\label{sec:cylacs}
An almost complex structure $J$ on $X$ is a bundle-endomorphism of $TX$ satisfying $J^2=-\id$. We say that an almost complex structure $J$ on a symplectic manifold $(X,\omega)$ is \emph{compatible with the symplectic form} if $\omega(\cdot,J\cdot)$ is a Riemannian metric on $X$, and \emph{tamed by the symplectic form} if $\omega(v,Jv) > 0$ whenever $v \neq 0$. It is well-known that the spaces of such almost complex structures are non-empty and contractible.

An almost complex structure on the symplectisation $(\R \times Y,d(e^t\lambda))$ of $(Y,\lambda)$ is said to be \emph{cylindrical} if it is compatible with the symplectic form, invariant under translations of the $t$-coordinate, satisfies $J\partial_t=R$, and preserves the contact-planes $\ker \lambda \subset TY$.

\subsubsection{The almost complex structure induced  by a metric}
\label{sec:acsmetric}
Assume that we are given a Lagrangian immersion
\[\iota \co L \to (X,\omega)\]
with transverse double-points. Here we describe how to construct a compatible almost complex structure in a neighbourhood of $\iota(L)$ induced by the choice of a Riemannian metric on $L$.

By Weinstein's Lagrangian neighbourhood theorem, the immersion $\iota$ can be extended to a symplectic immersion
\[\widetilde{\iota} \co (D^*L,d\theta_L) \to (X,\omega)\]
of the co-disc bundle $(D^*L,d\theta_L) \subset (T^* L, d\theta_L)$, with fibres of sufficiently small radius. We require that $\widetilde{\iota}$ restricts to $\iota$ along the zero-section.

Constructing $\widetilde{\iota}$ with some care, we can moreover make it satisfy the following property. Let $\mathcal{Q} \subset \Pi_{\mathrm{Lag}}(L)$ be the set of double-points. We require there to be disjoint neighbourhoods
\begin{gather*}
p_i \in U_{q,i} \subset L, \:\: i=1,2,\\
\{ p_1,p_2\} = \Pi_{\mathrm{Lag}}^{-1}(q),
\end{gather*}
for each $q \in \mathcal{Q}$, and diffeomorphisms $\varphi_{q,i} \co U_{q,i} \to D^n$ to the $n$-disc, such that
\begin{gather*}
(\widetilde{\iota}|_{D^*U_{q,2}} \circ \varphi_{q,2}^*)^{-1} \circ (\widetilde{\iota}|_{D^*U_{q,1}} \circ \varphi_{q,1}^*)\co D^*D^n \to D^*D^n,\\
(\mathbf{q},\mathbf{p}) \mapsto (-\mathbf{p},\mathbf{q}),
\end{gather*}
holds in some neighbourhood of the zero-section.

Let $g$ be the choice of a Riemannian metric on $L$ that coincides with the Euclidean metric in the above coordinates $\varphi_{q,i} \co U_{q,i} \to D^n$, where $U_{q,i}$ is a neighbourhood of the pre-image of a double-point. The metric $g$ determines a compatible almost complex structure $J_g$ on $(T^*L,d\theta_L)$ as follows.

First, the metric $g$ induces the Levi-Civita connection on the cotangent bundle $\pi \co T^*L \to L$
which, in turn, determines the horizontal subbundle $H_x \subset T_xT^*L$
for every $x \in T^*L$. The horizontal subbundle is a complement of the corresponding vertical subbundle $V_x := \ker (T\pi)_x.$ There are canonical identifications $H_x \simeq T_{\pi(x)}L$ and $V_x \simeq T^*_{\pi(x)}L$.

We define $J_g$ by requiring that $J_gH_x=V_x$ and that, for the horizontal vector identified with $h \in T_{\pi(x)}L$, the vertical vector $J_gh$ is identified with the covector $g(h,\cdot) \in T^*_{\pi(x)}L$. We refer to \cite[Remark 6.1]{DualityLeg} for an expression of this almost complex structure in local coordinates.

Finally, the almost complex structure $J_g$ can be pushed forward using $\widetilde{\iota}$ to a compatible almost complex structure in a neighbourhood of $\iota(L)$, as required. Observe that $J_g$ satisfies \ref{ra} since, by construction, there are holomorphic coordinates near each double-point in which the two intersecting sheets of $\iota(L)$ coincide with $\re(\C^n) \cup \im(\C^n)$.

\subsubsection{Adjusted triples $(f,g,J)$}
\label{sec:acs}
Again, we consider a Lagrangian immersion $\iota(L) \subset X$ as above. Also, let $f \co L \to \R$ be a Morse function, $g$ a Riemannian metric on $L$, and $J$ an almost complex structure on $X$. We say that the triple $(f,g,J_P)$ is \emph{adjusted to $\iota(L)$}  (see \cite[Section 6.3.1]{ContHomP}) given that the following holds.
\begin{itemize}
\item In a neighbourhood of $\iota(L) \subset X$, the compatible almost complex structure $J_P$ is induced by the metric $g$ on $L$ as above. In particular, we require each pre-image of a double-point to have a neighbourhood $U_{q,i} \subset L$ with coordinates $\varphi_{q,i} \co U_{q,i} \to D^n$ in which $g$ is the Euclidean metric.
\item The function $f$ is real-analytic and without critical points in the above coordinates $\varphi_{q,i} \co U_{q,i} \to D^n$ near the pre-image of a double-point.
\item The eigenvalues of the Hessian of $f$ (with respect to the metric $g$) at a critical point all have the same absolute values.
\end{itemize}

By \cite[Lemma 6.5]{ContHomP} it follows that there are triples $(f,g,J_P)$ adjusted to $\Pi_{\mathrm{Lag}}(\Lambda)$ for which $(f,g)$ is a Morse-Smale pair, and such that $J_P$ is regular. The importance of adjusted triples comes from the fact that \cite[Theorem 3.6]{ContHomP} applies to them. This establishes an identification of $J_P$-holomorphic discs on the two-copy of $\Pi_{\mathrm{Lag}}(\Lambda)$ with (generalised) $J_P$-holomorphic discs on $\Pi_{\mathrm{Lag}}(\Lambda)$, in case when the two-copy is obtained by pushing $\Pi_{\mathrm{Lag}}(\Lambda)$ off itself using $df \subset D^*\Lambda$. We refer to Section \ref{sec:twocopy} for more details.

\section{Background on Legendrian contact homology}
\label{sec:background-lch}

In the following, we assume that $\Lambda \subset P \times \R$ is a closed Legendrian submanifold which is \emph{chord generic}, that is, $\Pi_{\mathrm{Lag}}(\Lambda)$ is a generic immersion whose self-intersections thus consist of transverse double-points. Observe that this always can be assumed to hold after an arbitrarily $C^1$-small Legendrian isotopy. In particular, it follows that the set $\mathcal{Q}(\Lambda)$ of Reeb chords on $\Lambda$ is finite.

Recall that we require there to exist compatible almost complex structures $J_P$ on $P$ for which $(P,d\theta,J_P)$ has finite geometry at infinity, as defined in \cite[Definition 2.1]{ContHomP}. Also, we will only consider almost complex structures on $P$ of this kind. Given that one such almost complex structure exists, this condition is however not too restrictive. Namely, any deformation of such a compatible almost complex structure still satisfies this property, given that the deformation is supported inside a compact subset of $P$.
 
\subsection{The grading}
\label{sec:grading}

To each Reeb chord $c$ on $\Lambda$ we associate a grading
\[ | c | = CZ(\Gamma_c)-1, \]
where $\Gamma_c$ is a path of Lagrangian tangent-planes in $\C^n$ associated to $c$, and where $CZ(\Gamma_c)$ denotes the Conley-Zehnder index of this path as defined in \cite{ContHomR}. The path $\Gamma_c$ is obtained as follows.

One says that a Reeb chord is \emph{pure} if its endpoints lie on the same component of $\Lambda$ and, otherwise, one says that it is \emph{mixed}.

In the case when $c$ is a pure Reeb chord we let $\Gamma_c$ be the tangent-planes of $\Pi_{\mathrm{Lag}}(\Lambda)$ along the choice of a \emph{capping path} $\gamma_c$ of $c$, which is a continuous path on $\Lambda$ with starting-point (respectively end-point) at the end-point (respectively starting-point) of $c$. Furthermore, we assume that $\Pi_{\mathrm{Lag}}(\gamma_c)$ is null-homologous in $P$ and we choose a symplectic trivialisation of $TP$ along $\Pi_{\mathrm{Lag}}(\gamma_c)$ induced by a chain bounding $\Pi_{\mathrm{Lag}}(\gamma_c)$.

This construction provides a well-defined grading of the pure Reeb chords modulo the Maslov number of $\Lambda$ and twice the first Chern number of $P$. For simplicity, we will in the following assume both to be zero.

For a mixed Reeb chord there is obviously no capping path in the above sense. Instead, we proceed as follows. Fix points $p$, $q$ on two different components of $\Lambda$, together with a path $\gamma$ in $P$ connecting $p$ and $q$. Also, fix the choice of a path of Lagrangian planes in $TP$ along $\gamma$ starting at $T_p \Pi_{\mathrm{Lag}}(\Lambda)$ and ending at $T_q \Pi_{\mathrm{Lag}}(\Lambda)$.

For each mixed Reeb chord $c$ from the component containing $p$ to the component containing $q$, we construct a capping path by choosing a path on $\Lambda$ from the end-point of $c$ to $p$, followed by $\gamma$, and finally followed by a choice of path on $\Lambda$ from $q$ to the starting-point of $c$. Joining the corresponding paths of Lagrangian tangent-planes of $\Pi_{\mathrm{Lag}}(\Lambda)$ with the above choices of Lagrangian planes along $\gamma$, we get the desired path $\Gamma_c$.

It should be pointed out that the grading of the mixed Reeb chords starting and ending at two fixed components of $\Lambda$ depend on the above choice of a curve $\gamma$ together with the path of Lagrangian planes along this curve. However, given two Reeb chords both whose starting-points and end-points are located on $\Lambda_1$ and $\Lambda_2$, respectively, the \emph{difference in grading} between these two chords is independent of these choices.

\subsection{The relevant moduli spaces}
The differential in Legendrian contact homology is defined by a count of certain pseudo-holomorphic discs. We begin with the definitions of the moduli spaces that contain these discs.

\subsubsection{The moduli spaces of pseudo-holomorphic polygons}

\label{sec:moduli-polygons}
Fix a compatible almost complex structure $J_P$ on $P$. Given double-points $a,b_1,\hdots,b_m$ of $\Pi_{\mathrm{Lag}}(\Lambda)$, and writing $\mathbf{b}:=b_1\cdot \hdots \cdot b_m$, we let
\[\mathcal{M}_{a;\mathbf{b}}(\Pi_{\mathrm{Lag}}(\Lambda);J_P)\]
denote the moduli space of continuous maps
\[u \co (D^2,\partial D^2) \to (P,\Pi_{\mathrm{Lag}}(\Lambda)),\]
which are smooth in the interior, where they moreover satisfy
\[ \overline{\partial}_{J_P}(u):=du+J_P du \circ i =0.\]
Furthermore, we require there to be $m+1$ distinct boundary points $p_0, \hdots, p_m$, appearing in this cyclic order relative to the boundary orientation, such that the following holds.
\begin{itemize}
\item The map $u|_{\partial D^2 \setminus \{p_i\}}$ has a continuous lift to $\Lambda$ under $\Pi_{\mathrm{Lag}}$.
\item $u$ maps $p_0$ to $a$ and the $z$-coordinate of the above lift of $u|_{\partial D^2 \setminus \{p_i\}}$ is required to make a positive jump when traversing $p_0$ in positive direction according to the boundary orientation.
\item $u$ maps $p_i$ to $b_i$ for $i>0$ and the $z$-coordinate of the above lift of $u|_{\partial D^2 \setminus \{p_i\}}$ is required to make a negative jump when traversing $p_i$ in positive direction according to the boundary orientation.
\end{itemize}
Finally, two solutions are identified if they differ by a biholomorphism of the domain. We refer to $p_0$ as a \emph{positive boundary puncture mapping to $a$}, and to $p_i$ with $i>0$ as a \emph{negative boundary puncture mapping to $b_i$}. We will refer to the above above discs as \emph{pseudo-holomorphic polygons}.

\subsubsection{The moduli spaces of pseudo-holomorphic discs with strip-like ends in the symplectisation}
\label{sec:moduli-discs}
Here we assume that $V \subset (\R \times (P \times \R),d(e^t(dz+\theta)))$ is an exact Lagrangian cobordism from $\Lambda_-$ to $\Lambda_+$.

Fix a compatible almost complex structure $J$ on $\R \times (P \times \R)$ which we require to be cylindrical outside of some set of the form $[-N,N] \times (P \times \R)$. Given Reeb chords $a \in \mathcal{Q}(\Lambda_+)$ and $b_1,\hdots,b_m \in \mathcal{Q}(\Lambda_-)$, and writing $\mathbf{b}:=b_1\cdot \hdots \cdot b_m$,  we let
\[\mathcal{M}_{a;\mathbf{b}}(V;J)\]
denote the moduli space of continuous maps
\[\widetilde{u} =(\alpha,u) \co (\dot{D}^2,\partial \dot{D}^2) \to (\R \times (P \times \R),V)\]
which are smooth in the interior, where they moreover satisfy
\[ \overline{\partial}_J(\widetilde{u}):=d\widetilde{u}+J d\widetilde{u} \circ i =0.\]
Here $\dot{D}^2=D^2\setminus \{p_0,\hdots,p_m\}$ for a fixed choice of $m+1$ distinct points $p_0,\hdots,p_m \in \partial D^2$, appearing in this cyclic order relative to the boundary orientation, and $i$ is some conformal structure on $\dot{D}^2$.

We moreover require that $\widetilde{u}$ has the following asymptotical behaviour. Let $\gamma_c(t) \co [0,\ell(c)] \to P \times \R$ be the parametrisation of a Reeb chord $c$ on $\Lambda$ for which $\gamma_c'(t)=\partial_z$. Also, for each Reeb chord on $\Lambda$, we fix a contact-form preserving identification of a neighbourhood of the Reeb chord with $(D^*D^n \times \R,dz+\theta_{D^n})$.
\begin{itemize}
\item Let $s+it$ be coordinates on $\dot{D}^2$ induced by a biholomorphic identification of $D^2$ with the strip $\{s+it; \: 0 \le t \le \ell(a)\} \subset \C$ under which $p_0$ corresponds to $s=+ \infty$. We require there to be $s_0 \in \R$ and $\lambda >0$ such that
\[\| (\alpha,u)(s+it)-(s+s_0,\gamma_a(t)) \| \le e^{-\lambda|s|}\]
holds in the above coordinates, for all $s > 0$ sufficiently large.
\item Let $s+it$ be coordinates on $\dot{D}^2$ induced by a biholomorphic identification of $D^2$ with the strip $\{s+it; \: 0 \le t \le \ell(b_i)\} \subset \C$ under which $p_i$ corresponds to $s=- \infty$. We require there to be $s_0 \in \R$ and $\lambda >0$ such that
\[\| (\alpha,u)(s+it)-(s+s_0,\gamma_{b_i}(t)) \| \le e^{-\lambda|s|}\]
holds in the above coordinates, for all $s < 0$ sufficiently small.
\end{itemize}
Finally, two solutions are identified if they differ by a biholomorphism of the domain. We refer to $p_0$ as a \emph{positive boundary puncture asymptotic to $a$}, and to $p_i$ with $i>0$ as a \emph{negative boundary puncture asymptotic to $b_i$}.

In the definition of the Chekanov-Eliashberg algebra we will consider the cylindrical cobordism $V=\R \times \Lambda$ and an almost complex structure $J$ which is cylindrical on all of the symplectisation. Observe that, in this case, the above moduli spaces carry a natural $\R$-action induced by translation of the $t$-coordinate.

\subsubsection{Energies}
\label{sec:energy}

We will use the following notions of energies for the pseudo-holomorphic discs of the two kinds considered above.

For a disc $u \co D^2 \to P$, we define its symplectic area, also called $d\theta$-energy, by
\[E_{d\theta}(u) := \int_u d\theta. \]

Consider a $J$-holomorphic disc $\widetilde{u} \co D^2 \to \R \times (P \times \R)$ having boundary on $\R \times \Lambda$, and where $J$ is cylindrical. Following \cite[Section 5.3]{CompSFT}, we define its $d\theta$ and $\lambda$-energy by
\begin{eqnarray*}
& & E_{d\theta}(\widetilde{u}) := \int_{\widetilde{u}} d\theta, \\
& & E_\lambda(\widetilde{u}) := \sup_{\rho \in \mathcal{C}}\int_{\widetilde{u}} \rho(t) dt \wedge (dz+\theta),
\end{eqnarray*}
respectively, where $\mathcal{C}$ is the set of smooth functions $\rho \co \R \to \R_{\ge 0}$ having compact support and satisfying $\int_\R \rho(t) dt=1$.

To each Reeb chord $c$, we assign the action
\[ \ell(c) := \int_c (dz+\theta).\]
Given $u \in \mathcal{M}_{a;\mathbf{b}}(\R \times \Lambda;J) $ and $\widetilde{u} \in \mathcal{M}_{a;\mathbf{b}}(\Pi_{\mathrm{Lag}}(\Lambda);J_P)$, where $J$ again is assumed to be cylindrical, a calculation similar to \cite[Lemma B.3]{RationalSFT} and \cite[Lemma 5.16]{CompSFT} yields
\begin{eqnarray*}
& & 0 \leq E_{d\theta}(u)=E_{d\theta}(\widetilde{u})=\ell(a) -(\ell(b_1)+\hdots+\ell(b_m)), \\
& & 0 < E_\lambda(\widetilde{u})=\ell(a).
\end{eqnarray*}

Observe that the $d\theta$-energy of a pseudo-holomorphic disc in $P$ vanishes if and only if the disc is constant, while the $d\theta$-energy of a pseudo-holomorphic disc in the symplectisation vanishes if and only if it is contained entirely in a trivial strip $\R \times c$ over a Reeb chord $c$.

\subsubsection{Dimension formulae}
\label{sec:dim}
We call an almost complex structure \emph{regular} if the appropriate moduli spaces are transversely cut out, and hence are smooth finite-dimensional manifolds.  In this case, it follows from the calculation of the Fredholm index in \cite[Section 6]{ContHomR} (also see \cite[Proposition 2.3]{ContHomP}), that the dimensions of the above moduli spaces are given by
\begin{gather*}
\dim \mathcal{M}_{a;\mathbf{b}}(V;J)=|a|-|b_1|-\hdots-|b_m|,\\
\dim \mathcal{M}_{a;\mathbf{b}}(\Pi_{\mathrm{Lag}}(\Lambda);J_P)=|a|-|b_1|-\hdots-|b_m|-1.
\end{gather*}
Recall that we here assume that the Maslov classes of $V$ and $\Lambda$, as well as the first Chern class of $P$, all vanish.

In particular,
observe that
\[\dim \mathcal{M}_{a;\mathbf{b}}(\R \times \Lambda;J) = \dim \mathcal{M}_{a;\mathbf{b}}(\Pi_{\mathrm{Lag}}(\Lambda);J_P) +1.\]
In the case when $J$ is cylindrical, the extra degree of freedom in the moduli space of $J$-holomorphic discs on the left-hand side should be thought of as coming from the translations of the $t$-coordinate.

Since the pseudo-holomorphic discs in $\mathcal{M}_{a;\mathbf{b}}(\Pi_{\mathrm{Lag}}(\Lambda);J_P)$ have only one positive puncture, any compatible almost complex structure on $(P,d\theta)$ satisfying \ref{ra} can be approximated by a regular almost complex structure \cite[Lemma 4.5]{ContHomP}.

Any solution $\widetilde{u} \in \mathcal{M}_{a;\mathbf{b}}(V;J)$ is necessarily always injective, i.e.~it has an interior-point $x \in D$ for which $u^{-1}(u(x))=\{x\}$ and $du(x) \neq  0$, as follows by its asymptotical properties. The standard transversality argument \cite[Chapter 3]{JholCurves} applies, showing that $J$ can be made regular after an arbitrarily small compactly supported perturbation.

Observe that more work is needed in order to find a regular $J$ that is also \emph{cylindrical}. See \cite{FredholmTheory} for the case without boundary. However, in our setting, we will use Theorem \ref{thm:lift} to deduce the regularity for the cylindrical lift $\widetilde{J}_P$ of a regular almost complex structure $J_P$ on $P$.

\subsection{The Chekanov-Eliashberg algebra}
Consider the unital graded $\Z_2$-algebra $\mathcal{A}_\bullet(\Lambda)$ freely generated by the Reeb chords on $\Lambda$ with grading determined by the Conley-Zehnder index as above. The Chekanov-Eliashberg algebra of $\Lambda$ is the chain complex $(\mathcal{A}_\bullet(\Lambda),\partial)$, where the differential $\partial$ is defined as follows.

Let $a$ be a Reeb chord generator of the algebra. The differential given in \cite{ContHomP} is defined by
\[ \partial(a) := \sum_{|a|-|\mathbf{b}|=1} |\mathcal{M}_{a;\mathbf{b}}(\Pi_{\mathrm{Lag}}(\Lambda);J_P)|\mathbf{b}\]
for some choice of regular compatible almost complex structure $J_P$ on $P$. Similarly, the differential given in \cite{RationalSFT} is defined by
\[ \partial(a) := \sum_{|a|-|\mathbf{b}|=1} |\mathcal{M}_{a;\mathbf{b}}(\R \times \Lambda;J)/\R|\mathbf{b}\]
for some choice of regular cylindrical almost complex structure $J$ on $\R \times (P \times \R)$. Observe that the dimension formula implies that the sum is taken over zero-dimensional moduli spaces. Together with the Gromov-Hofer compactness in \cite{ContHomP} and \cite{CompSFT}, it follows that the above counts make sense.

In both of the above cases the differential is extended to the whole algebra via the Leibniz rule
\[ \partial(ab)=\partial(a)b+a\partial(b),\]
and it follows that $\partial$ is of degree $-1$. Moreover, the above formula for the $d\theta$-energy implies that the differential is action-decreasing, and hence that there is an induced filtration of the complex $\mathcal{A}_\bullet$ induced by the action.

We now refer to the invariance results for the above two versions of Legendrian contact homology, that is, the following holds for either of the above two definitions of the boundary operator $\partial$.
\begin{theorem}[Theorem 1.1 in \cite{ContHomP}, \cite{RationalSFT}]
For a closed Legendrian submanifold $\Lambda \subset P \times \R$ of the contactisation of a Liouville domain, it is the case that
\begin{itemize}
\item $\partial^2=0$.
\item The homotopy type of $(\mathcal{A}_\bullet(\Lambda),\partial)$ is independent of the choice of a regular compatible almost complex structure, and invariant under Legendrian isotopy.
\end{itemize}
\end{theorem}

\subsection{Linearised Legendrian contact homology}
\label{sec:linearised}

The results in \cite[Lemma 3.15]{RationalSFT} and \cite[Section 4]{RationalSFT} show that an exact Lagrangian cobordism $V$ from $\Lambda_-$ to $\Lambda_+$ induces a unital DGA-morphism
\[ \Phi_V \co (\mathcal{A}(\Lambda_+),\partial_+) \to (\mathcal{A}(\Lambda_-),\partial_-).\]
It is defined by counting rigid $J$-holomorphic discs in $\R \times (P \times \R)$ having boundary on $L$ and boundary-punctures asymptotic to Reeb chords. Here we require $J$ to coincide with the cylindrical almost complex structures $J_{\pm\infty}$ in the sets $\{t \ge N\}$ and $\{ t \le -N \}$, respectively, where $N > 0$ is sufficiently large, and $J_{\pm\infty}$ are used in the definitions of the above Chekanov-Eliashberg algebras. Again, we assume that the Maslov class of $V$ vanishes.

In particular, as is also shown in \cite{RationalSFT2}, an exact Lagrangian filling $L$ of $\Lambda$ together with an appropriate choice of almost complex structure induces a unital DGA-morphism
\[\epsilon_L \co (\mathcal{A}(\Lambda),\partial) \to (\Z_2,0),\]
where the right-hand side is the trivial DGA. In general, a unital DGA-morphism to $(\Z_2,0)$ is called an \emph{augmentation}.

Given an augmentation $\epsilon$ of a semi-free DGA $(\mathcal{A}_\bullet,\partial)$, one can construct the following chain complex. Define an algebra automorphism $\Psi^\epsilon$ of $\mathcal{A}_\bullet$ by prescribing $\Psi^\epsilon(a)=a+\epsilon(a)$ for each generator $a$. It follows that the constant part of $\Psi^\epsilon\circ\partial \circ (\Psi^\epsilon)^{-1}$ vanishes and, consequently, its linear part
\[\partial_\epsilon:=(\Psi^\epsilon\circ\partial \circ (\Psi^\epsilon)^{-1})^1\]
is itself a differential on the graded $\Z_2$-vector space $\mathcal{A}^1_\bullet$ spanned by the generators of $\mathcal{A}_\bullet$. We call the chain complex $(\mathcal{A}^1_\bullet,\partial_\epsilon)$ the \emph{linearisation of the DGA induced by $\epsilon$}.

Obviously, an augmentation of $(\mathcal{A}_\bullet,\partial)$ pulls back under a unital DGA-morphism $\Phi \co (\mathcal{A}'_\bullet,\partial') \to (\mathcal{A}_\bullet,\partial)$ to an augmentation $\epsilon':=\epsilon\circ \Phi$ of $(\mathcal{A}'_\bullet,\partial')$. Moreover, the map $\Psi^\epsilon \circ \Phi \circ  (\Psi^{\epsilon'})^{-1}$ has vanishing constant part. We denote its linear part by
\[\Phi_\epsilon:=(\Psi^\epsilon \circ \Phi \circ  (\Psi^{\epsilon'})^{-1})^1 \co (\mathcal{A}'^1_\bullet,\partial'_{\epsilon'}) \to (\mathcal{A}^1_\bullet,\partial_\epsilon).\]
This can be seen to be a chain-map between the corresponding linearisations.

If the Chekanov-Eliashberg algebra of a Legendrian $\Lambda$ has an augmentation $\epsilon$, we denote the corresponding linearisation by
\[(CL_\bullet(\Lambda),\partial_\epsilon):=(\mathcal{A}_\bullet(\Lambda)^1=\Z_2\langle \mathcal{Q}(\Lambda) \rangle,\partial_\epsilon),\]
also called its \emph{linearised Legendrian contact homology complex}. We use $HCL_\bullet(\Lambda;\epsilon)$ to denote its homology. We will also use $(CL^\bullet(\Lambda),d_\epsilon)$ to denote the associated co-complex, and $HCL^\bullet(\Lambda;\epsilon)$ to denote the corresponding cohomology. 

The homotopy type of this complex does indeed depend on the choice of augmentation. However, \cite[Theorem 5.2]{DiffAlg} shows that the set of all isomorphism classes of linearised homologies is a Legendrian isotopy invariant. This proof however depends on the invariance result for the version of the Chekanov-Eliashberg algebra as defined in \cite{DiffAlg} and \cite{ContHomP}, which gives more than just invariance of the homotopy type; it establishes that the \emph{stable tame isomorphism type} of the Chekanov-Eliashberg algebra is a Legendrian isotopy invariant.

\begin{remark}
\label{rem:stabletame}
A priori it is not clear from the invariance proof in \cite{RationalSFT} that the stable tame isomorphism type of the Chekanov-Eliashberg algebra defined using the symplectisation is an invariant. Of course, Theorem \ref{thm:lift} is one way to establish this result in the latter setting as well. In particular, the invariance result \cite[Theorem 5.2]{DiffAlg} for the linearised Legendrian contact homology applies here as well.
\end{remark}

We now recall the invariance property for the linearised contact homology induced by a filling.
\begin{theorem}[Theorems 1.1 and 2.1 in \cite{RationalSFT2}]
\label{thm:lininv}
Let $L \subset \R \times (P \times \R)$ be a filling of $\Lambda$. The homotopy type of $(CL^\bullet(\Lambda),d_{\epsilon_L})$ is independent of the choice of an almost complex structure. Suppose that $L'$ is a filling of $\Lambda'$ that is isotopic to $L$ by a Hamiltonian isotopy supported in a set of the form $\R \times K$, where $K \subset P \times \R$ is compact. It follows that there is a homotopy equivalence
\[(CL^\bullet(\Lambda),d_{\epsilon_L}) \simeq (CL^\bullet(\Lambda'),d_{\epsilon_{L'}}).\]
\end{theorem}
In the setting considered here, this invariance result can also be seen to follow from techniques similar to the ones used in the proof of Theorem \ref{thm:twocopy}, together with the invariance of wrapped Floer homology under Hamiltonian isotopies (see Theorem \ref{thm:inv}).

\section{Background on wrapped Floer homology}
\label{sec:background-floer}

We will now give an outline of wrapped Floer homology as defined in \cite{RationalSFT2}. Wrapped Floer homology is a Hamiltonian isotopy invariant of pairs of exact Lagrangian fillings inside an exact symplectic manifold. Here we will only consider the case when the ambient symplectic manifold is the symplectisation of a contactisation, even though the theory is defined in more generality.

In the following, we thus let $L$ and $L'$ be exact Lagrangian fillings of $\Lambda$ and $\Lambda'$, respectively, inside the symplectisation of a contactisation
\[(\R \times (P \times \R),d(e^t\lambda)),\:\: \lambda:=dz+\theta,\]
where $P$ is $2n$-dimensional. To simplify the definition of the grading, we assume that the Maslov class vanishes for both $L$ and $L'$, and that $P$ has vanishing first Chern class.

\subsection{Definition of the wrapped Floer homology complex}
After a generic Hamiltonian perturbation of $L'$ we may assume that $L$ and $L'$ intersect transversely in finitely many double-points, and that the Legendrian link $\Lambda \cup \Lambda'$ is embedded and chord-generic.

\subsubsection{The graded vector space}
Use $c_1,\hdots,c_l$ to denote the Reeb chords starting on $\Lambda$ and ending on $\Lambda'$, and $x_1,\hdots,x_k$ to denote the double-points of $L \cup L'$. We define the vector spaces
\begin{eqnarray*}
CF^0(L,L') & := & \Z_2\langle x_1,\hdots,x_k \rangle,\\
CF^\infty(L,L') & := & \Z_2 \langle c_1, \hdots, c_l \rangle,\\
CF(L,L') & := & CF^0(L,L') \oplus CF^\infty(L,L'),
\end{eqnarray*}
which we endow with the following grading.

Assume that both $L$ and $L'$ are cylindrical in the set $\{ t \ge N\}$  and consider a function $\sigma(t) \co \R \to \R_{\ge 0}$ satisfying $\sigma'(t) \ge 0$, $\sigma(t)=0$ for $t \le N$, and $\sigma(t) =1$ for all sufficiently big $t > 0$. There is a Hamiltonian $H \co \R \times (P \times \R) \to \R$ depending only on the $t$-coordinate, and whose Hamiltonian vector field is given by
\[X_H=-\sigma(t)\partial_z.\]
Its time-$s$ flow $\phi^s_H$ has the property that, for generic $\sigma(t)$ and $s > 0$ sufficiently large, the double-points of $L \cup \phi^s_H(L')$ appearing in $\{ t \ge N \}$ are transverse double-points which naturally are in bijective correspondence with the Reeb chords starting on $\Lambda$ and ending on $\Lambda'$. To see this, observe that $\phi^s_H$ wraps $L'$ in the negative Reeb-direction.

We fix a Legendrian lift of the exact Lagrangian immersion $L \cup \phi^s_H(L')$ to the contactisation of the symplectisation, which moreover has the property that all Reeb chords start on the lift of $L$. Using the constructions in Section \ref{sec:grading}, we can associate paths $\Gamma$ of Lagrangian tangent-planes in the symplectisation to each generator (considered as a mixed Reeb chord on the lift to the contactisation of the symplectisation). The gradings of the generators are then defined to be
\[|x_i|:=CZ(\Gamma_{x_i}), \:\: |c_i|:=CZ(\Gamma_{c_i}).\]

\begin{remark}
Observe that for a Reeb-chord generator $c_i$, this grading differs from the grading $|c_i|_{LCH}$ in Section \ref{sec:grading} obtained when considering $c_i$ as a generator of the Chekanov-Eliashberg algebra of $\Lambda \cup \Lambda'$. More precisely, the capping paths in $P$ used for the grading $|c_i|_{LCH}$ can be lifted to capping paths in the symplectisation, and the corresponding gradings are related by $|c_i|=|c_i|_{LCH}+2$.
\end{remark}

\begin{remark}
Recall that the above grading is not canonical; the different choices of paths of Lagrangian planes along the curve $\gamma$ in the construction given in Section \ref{sec:grading} may induce a global shift in grading of $CF_\bullet$. However, in the case when $L'$ is a sufficiently $C^1$-small perturbation of $L$, which moreover is assumed to be connected, a canonical grading is obtained as follows; we choose both the path $\gamma$, and the path of Lagrangian planes along $\gamma$, to be sufficiently close to constant paths.
\end{remark}

\subsubsection{The differential}
Under the decomposition
\[CF_\bullet(L,L')=CF_\bullet^0(L,L') \oplus CF_\bullet^\infty(L,L')\]
the differential will be of the form
\[ \partial := \begin{pmatrix}
\partial_0 & 0 \\
\delta & \partial_\infty
\end{pmatrix},
\]
where the entries in the matrix are to be defined below. In other words, we will construct the wrapped Floer homology complex as the mapping cone
\[(CF(L,L'),\partial):=\Cone(\delta)\]
of a chain-map $\delta$.

We fix a compatible almost complex structure $J$ on $\R \times (P \times \R)$ which coincides with the cylindrical almost complex structure $J_\infty$ in the set $\{ t \ge N\}$ for some sufficiently big $N > 0$. In the following we assume that $J$ is regular for the below spaces of pseudo-holomorphic discs. To see how this can be achieved, we refer to Section \ref{sec:reg} below.

We will consider pseudo-holomorphic discs having boundary on $L \cup L'$ and boundary-punctures of which some are asymptotic to Reeb chords on the Legendrian end, and some are mapped to double-points of $L \cup L'$. We will require that these pseudo-holomorphic discs satisfy the conditions in Section \ref{sec:moduli-discs} outside of some compact set, and that they satisfy the conditions in Section \ref{sec:moduli-polygons} in some neighbourhood of the punctures which are mapped to double-points.

Recall that we have fixed a Legendrian lift of $L \cup L'$ with the property that every Reeb chord starts on $L$ and ends on $L'$, and that this choice induces a notion of positivity and negativity for the punctures of the latter kind.

\subsubsection{The sub-complex $CF^\infty$}
The two fillings $L$ and $L'$, together with the almost complex structure $J$, induce augmentations $\epsilon_L$ and $\epsilon_{L'}$ of the Chekanov-Eliashberg algebra of $\Lambda$ and $\Lambda'$, respectively. Here we let both DGAs be defined using the cylindrical almost complex structure $J_\infty$.

There is an induced augmentation $\epsilon_{L \cup L'}$ of the Chekanov-Eliashberg algebra of the Legendrian link $\Lambda \cup \Lambda'$ which vanishes on generators corresponding to chords between $\Lambda$ and $\Lambda'$, and which takes the value $\epsilon_L$ and $\epsilon_{L'}$ on generators corresponding to chords on $\Lambda$ and $\Lambda'$, respectively.

The Reeb chords on $\Lambda \cup \Lambda'$ starting on $\Lambda$ and ending on $\Lambda'$ span a sub-complex of the linearised Legendrian contact homology complex of the link $\Lambda \cup \Lambda'$. We will use
\[(CL_\bullet(\Lambda, \Lambda'),\partial_{\epsilon_L,\epsilon_{L'}}) \subset (CL_\bullet(\Lambda \cup \Lambda'),\partial_{\epsilon_{L \cup L'}})\]
to denote this sub-complex. Also, we will let $(CL^\bullet(\Lambda, \Lambda'),d_{\epsilon_L,\epsilon_{L'}})$ denote the corresponding co-complex and we will use $HCL^\bullet(\Lambda, \Lambda';\epsilon_L,\epsilon_{L'})$ to denote the corresponding cohomology.

We define $\partial_\infty$ by making the identification
\[(CF_\bullet^\infty(L,L'),\partial_\infty):=(CL^{\bullet-2}(\Lambda, \Lambda'),d_{\epsilon_L,\epsilon_{L'}}),\]
where we recall the above grading conventions.

\subsubsection{The quotient-complex $CF^0$.}
We define the complex
\begin{gather*}
(CF^0(L,L'),\partial_0) := (\Z_2 \langle x_1, \hdots, x_k \rangle,\partial_0),\\
\partial_0(x_i)=\sum_{|x_j|-|x_i|=1}|\mathcal{M}_{x_j;x_i}(L \cup L';J)|x_j,
\end{gather*}
where $\mathcal{M}_{x_j;x_i}(L \cup L';J)$ denotes the moduli space of $J$-holomorphic polygons defined in Section \ref{sec:moduli-polygons}. Recall that $|x_j|-|x_i|-1$ is the dimension of this moduli space in the case when $J$ is regular. Since these polygons have only two punctures, we will sometimes refer to them as strips.

\begin{remark}
Recall that we are considering a Legendrian lift of $L \cup L'$ to the contactisation of the symplectisation for which all Reeb chords start on $L$ and end on $L'$. The complex $(CF^0(L,L'),\partial_0)$ is in fact nothing else than the linearised Legendrian contact cohomology complex of this lift. Note that, since there are no pure Reeb chords on this Legendrian lift, there is a canonical augmentation that maps every generator to zero.
\end{remark}

One can associate an action $\ell_{d(e^t{\lambda})}$ to every generator in $(CF^0(L,L'),\partial_0)$, by associating to it the action of the corresponding Reeb chord for the above choice of Legendrian lift of $L \cup L'$. It follows that a pseudo-holomorphic strip $u \in \mathcal{M}_{x_j;x_i}(L \cup L')$ has $d(e^t{\lambda})$-area given by
\[ 0 < E_{d(e^t{\lambda})}(u) = \ell_{d(e^t{\lambda})}(x_j)-\ell_{d(e^t{\lambda})}(x_i),\]
and that hence $\partial_0$ is action-increasing with respect to the action $\ell_{d(e^t{\lambda})}$.

\subsubsection{The chain map}
There is a chain map
\begin{gather*}
\delta \co (CF_\bullet^0(L,L'),\partial_0) \to (CL^{\bullet-1}(\Lambda, \Lambda'),d_{\epsilon_L,\epsilon_{L'}}),\\
\delta(x_i) := \sum_{|c_j|-|x_i|=1} |\mathcal{M}_{c_j;x_i}(L \cup L';J)| c_j,
\end{gather*}
where the moduli-space
\[\mathcal{M}_{c_j;x_i}(L \cup L';J)\]
consists of $J$-holomorphic discs inside $\R \times (P \times \R)$ having boundary on $L \cup L'$, a positive puncture asymptotic to the Reeb chord $c_j$ from $\Lambda$ to $\Lambda'$, and a negative puncture mapping to the double-point $x_i \in L \cap L'$.

Finally, observe that $|c_j|-|x_i|-1$ is the dimension of the above moduli space in the case when $J$ is regular.

\subsubsection{Finding regular almost complex structures}
\label{sec:reg}
The following method can be used to show the existence of almost complex structures that are regular for the above moduli spaces. Recall that $J$ is a compatible almost complex structure on $\R \times (P \times \R)$ which coincides with a cylindrical almost complex structure $J_\infty$ in the set $\{t \ge N \}$.

By the discussion in Section \ref{sec:dim}, after a compactly supported perturbation, we may assume that $J$ is regular for the discs in the definition of the augmentations. By Theorem \ref{thm:lift} it follows that, after choosing $J_\infty$ to be the cylindrical lift of a regular almost complex structure on $P$, we may assume that $J_\infty$ is regular for all discs in the definition of $\partial_\infty$.

It remains to show that the moduli spaces in the definition of $\partial_0$ and $\delta$ can be made transversely cut out. Observe that these discs have exactly one and two boundary punctures mapping to double-points, respectively. Assuming that $J$ satisfies \ref{ra} in a neighbourhood of the double-points $L \cap L'$, \cite[Lemma 4.5]{ContHomP} again applies to show that $J$ may be assumed to be regular for these moduli-spaces after a compactly supported perturbation.

\subsection{The transfer-map induced by an exact Lagrangian cobordism}
\label{sec:transfer}

Let $L$ and $L'$ be exact Lagrangian fillings of $\Lambda$ and $\Lambda'$ as above, and let $V \subset \R \times (P \times \R)$ be an exact Lagrangian cobordism from $\Lambda'$ to $\Lambda''$. We moreover assume that $L$ and $L'$ are cylindrical in the set $\{ t \ge -1 \}$.

Let $J$ denote the regular almost complex structure on $\R \times (P \times \R)$ defining $(CF_\bullet(L,L'),\partial)$, which we assume coincides with the cylindrical almost complex structure $J_\infty$ in the set $\{ t \ge -1\}$.

Assuming that $V$ is cylindrical in the set $\{ t \le 1 \}$, recall the definition of the concatenation
\[ L''_s:=L' \odot_s V, \:\: s \ge 0,\]
given in Section \ref{sec:defexact}, which is an exact Lagrangian filling of $\Lambda''$. We assume that $L$ and $L''$ intersect transversely. The double-points of $L \cup L''_s$ can be decomposed as
\[ L \cap L''_s = (L \cap L') \sqcup ((\R \times \Lambda) \cap V).\]

Fix an almost complex structure $J_V$ which coincides with $J_\infty$ on $\{t \le 1 \}$ and with the cylindrical almost complex structure $J''_\infty$ on some set of the form $\{t > N\}$. For each $s>0$, we consider the almost complex structure
\[J \odot_s J_V \]
on $\R \times (P \times \R)$ which coincides with $J$ in the set $\{ t \le s+1 \}$ and with $J_V(t-s,p,z)$ in the set $\{t \ge s+1\}$. We let $(CF_\bullet(L,L''_s),\partial''_s)$ be the induced wrapped Floer homology complex.

Let
\[CF^0_\bullet (\R \times \Lambda,V) \subset CF^0_\bullet (L,L' \odot_s V) \]
denote the subspace spanned by the double-points $(\R \times \Lambda) \cap V$. The above wrapped Floer homology complex is of the form
\begin{gather*}
CF_\bullet(L,L' \odot_s V)=CF^0_\bullet (L,L') \oplus CF^0_\bullet (\R \times \Lambda,V) \oplus CF^\infty_\bullet(L,L''_s),\\
\partial''_s=\begin{pmatrix} \partial_0 & \delta''_4 & 0 \\
\delta''_1 & \partial_V & 0 \\
\delta''_2 & \delta''_3 & \partial''_\infty
\end{pmatrix},
\end{gather*}
where $\partial_0$ is the differential of $CF^0_\bullet(L,L')$, given that $s >0$ is chosen sufficiently large.

To see the latter claim, observe that we can increase $s>0$ without changing the action of the generators in $CF^0_\bullet(L,L')$. A monotonicity argument for the $d(e^t\lambda)$-area of pseudo-holomorphic curves with boundary (see the proof of Lemma \ref{lem:mono}) can be used to give the following. The discs in the definition of $\partial_0$ are all contained in some set $\{ t \le A \}$ for $s>0$ sufficiently large, where $A$ is independent of $s$. Alternatively, this statement can be shown using a neck-stretching argument (see \cite[Section 3.4]{CompSFT}).
\begin{remark}
In the case when all generators of $CF^0_\bullet (\R \times \Lambda,V)$ have $\ell_{d(e^t \lambda)}$-action greater than $CF^0_\bullet (L,L')$, it immediately follows that $\delta''_4=0$.
\end{remark}

In \cite[Section 4.2.2]{RationalSFT2} the so called \emph{transfer map} is constructed which, for $s > 0$ sufficiently large, is a chain map of the form
\begin{gather*}
\Phi_V \co (CF_\bullet(L,L'),\partial) \to (CF_\bullet(L,L' \odot_s V),\partial''_s),\\
\Phi_V = \begin{pmatrix}
\id_{CF^0} & 0 \\
0 & \phi_0 \\
0 & \phi_\infty
\end{pmatrix},
\end{gather*}
relative the above decomposition.

To describe its components, we proceed as follows. Take generators
\begin{gather*}
c_i \in CF^\infty_\bullet (L,L'),\\
d_j \in CF^\infty_\bullet (L,L' \odot V),\\
x_j \in CF^0_\bullet (\R \times \Lambda,V),
\end{gather*}
and let $\mathbf{a}$ and $\mathbf{b}$ denote words of Reeb chords on $\Lambda'$ and $\Lambda$, respectively.

We use $\mathcal{M}_{d_j;\mathbf{a},c_i,\mathbf{b}}((\R \times \Lambda) \cup V;J_V)$ to denote the moduli space of $J_V$-holomorphic discs as defined in Section \ref{sec:moduli-discs} having boundary on $(\R \times \Lambda) \cup V$ and boundary punctures asymptotic to the prescribed Reeb chords.

Similarly, we define the moduli space $\mathcal{M}_{x_j;\mathbf{a},c_i,\mathbf{b}}((\R \times \Lambda) \cup V;J_V)$ consisting of $J_V$-holomorphic discs having boundary on $(\R \times \Lambda) \cup V$, a positive puncture mapping to $x_j$, and its negative punctures asymptotic to the prescribed Reeb chords. Here the Legendrian lift of of $(\R \times \Lambda) \cup V$ has been chosen so that all Reeb chords start on the lift of $\R \times \Lambda$. As in Section \ref{sec:moduli-polygons}, this induces the notions of positivity and negativity for a puncture mapping to a double-point.

\begin{remark}
\label{rem:nopos}
Observe that with this notion of positivity and negativity for the boundary punctures, it is not necessary for a $J_V$-holomorphic disc as above to posses a positive puncture. However, a strip without positive punctures must have a negative puncture asymptotic to a Reeb chord at $-\infty$ starting on $\Lambda'$ and ending on $\Lambda$.
\end{remark}

For a generator $c_i \in CF^\infty_\bullet(L,L')$, the components of $\Phi_V$ are given by the following counts of rigid $J_V$-holomorphic discs in the above moduli-spaces.
\begin{eqnarray*}
\phi_0(c_i) & := & \sum_{|x_j|-|c_i|=0 \atop |\mathbf{a}|=\mathbf{b}|=0}|\mathcal{M}_{x_j;\mathbf{a},c_i,\mathbf{b}}((\R \times \Lambda) \cup V;J_V)|\epsilon_{L'}(\mathbf{a})\epsilon_{L}(\mathbf{b})x_j,\\
\phi_\infty(c_i) & := & \sum_{|d_j|-|c_i|=0 \atop |\mathbf{a}|=\mathbf{b}|=0}|\mathcal{M}_{d_j;\mathbf{a},c_i,\mathbf{b}}((\R \times \Lambda) \cup V;J_V)|\epsilon_{L'}(\mathbf{a})\epsilon_{L}(\mathbf{b})d_j.
\end{eqnarray*}
In the case when $J_V$ is regular, the dimensions of the above moduli spaces are given by $|x_j|-|c_i|-|\mathbf{a}|-|\mathbf{b}|$ and $|d_j|-|c_i|-|\mathbf{a}|-|\mathbf{b}|$, respectively.

A neck-stretching argument can be used to show the following.
\begin{proposition}
\label{prop:concat}
Let $L$ and $L'$ be exact Lagrangian fillings. Given exact Lagrangian cobordisms $V$ and $W$, where the positive end of $V$ equals the negative end of $W$, it follows that
\[ \Phi_{V \odot W} = \Phi_W \circ \Phi_V \co (CF_\bullet(L,L'),\partial) \to (CF_\bullet(L,L' \odot V \odot W ),\partial''),\]
given that all involved almost complex structures have been appropriately chosen.
\end{proposition}

\subsection{Invariance under Hamiltonian isotopy}

The Wrapped Floer homology complex satisfies the following invariance property.

\begin{theorem}{\cite[Theorem 4.10]{RationalSFT2}}
\label{thm:inv}
Let $L, L', L'' \subset \R \times (P \times \R)$ be exact Lagrangian fillings. Suppose that $L''=\phi^1_{H_s}(L')$ is isotopic to $L'$ by a Hamiltonian isotopy supported in $\R \times K$, where $K$ is compact. There is a homotopy equivalence
\[\Phi \co (CF_\bullet(L,L'),\partial) \to (CF_\bullet(L,L''),\partial'').\]
\end{theorem}

Since the proof of Theorem \ref{thm:twocopy} is based on some constructions used in the proof of this invariance theorem, we formulate its main ingredients below. The core of the argument is as follows. Subsequent applications of Lemma \ref{lem:hamiltonian} implies the following standard fact. Let $L'$ and $L''$ be fillings of $\Lambda'$ and $\Lambda''$, respectively. Given a Hamiltonian isotopy $\phi^s_{H_s}$ as above, there are exact Lagrangian cobordisms $U$, $V$, and $W$ satisfying the following.
\begin{itemize}
\item $L' \odot U$ is isotopic to $L''=\phi^1_{H_s}(L')$ by a compactly supported Hamiltonian isotopy.
\item $U \odot V$ is isotopic to $\R \times \Lambda'$ by a compactly supported Hamiltonian isotopy.
\item $V \odot W$ is isotopic to $\R \times \Lambda''$ by a compactly supported Hamiltonian isotopy.
\end{itemize}
Thea idea is now to use transfer-maps induced by $U$, $V$, and $W$, together with an invariance result for compactly supported Hamiltonian isotopies. To that end, the following two propositions are needed.

\begin{proposition}{\cite[Section 4.2.1]{RationalSFT2}}
\label{prop:path}
Let $L$ and $L'$ be exact Lagrangian fillings. Given paths $\phi^s_{H_s}(L')$ of fillings and $J_s$ of almost complex structures, both which are fixed outside of some compact set, there is an induced homotopy equivalence
\[\Phi_{H_s,J_s} \co (CF_\bullet(L,L'),\partial) \to (CF_\bullet(L,\phi^1_{H_s}(L')),\partial'').\]
Here the former complex is defined using $J_0$ and the latter is defined using $J_1$. Furthermore, the restriction
\[\Phi_{H_s,J_s} |_{CF^\infty } \co (CF^\infty_\bullet(L,L'),\partial_\infty) \to (CF^\infty_\bullet(L,\phi^1_{H_s}(L')),\partial''_\infty)\]
is an isomorphism of complexes. In the case when there are no births or deaths of double-points during the Hamiltonian isotopy, it follows that $\Phi_{H_s,J_s}$ is an isomorphism of complexes as well.
\end{proposition}

\begin{proposition}{\cite[Section 4.2.3]{RationalSFT2}}
Let $L$ and $L'$ be exact Lagrangian fillings. Assume that we are given paths $\phi^s_{H_s}(V)$ of exact Lagrangian cobordisms, and $J_s:=J \odot_\sigma J_{V,s}$ of almost complex structures, both which are fixed outside of some compact set. Write $W:=\phi^1_{H_s}(V)$ and let the transfer maps $\Phi_V$ and $\Phi_W$ be defined using $J_{V,0}$ and $J_{V,1}$, respectively. The diagram
\[\xymatrix{
(CF_\bullet(L,L'),\partial) \ar@{=}[d] \ar[r]^-{\Phi_V} & (CF_\bullet(L,L' \odot_\sigma V),\partial'') \ar[d]^{\Phi_{H_s,J_s}}\\
(CF_\bullet(L,L'),\partial) \ar[r]^-{\Phi_W} & (CF_\bullet(L,L' \odot_\sigma W),\partial''')
}
\]
commutes up to homotopy, given that $\sigma>0$ is sufficiently large and that the positive end of $L'$ agrees with the negative ends of $V$ and $W$.
\end{proposition}

\subsection{The transfer map induced by the negative Reeb-flow}
\label{sec:Reebflow}
Let $L$ and $L'$ be fillings of $\Lambda$ and $\Lambda'$, respectively, as above. We moreover assume that both fillings are cylindrical in the set $\{ t \ge -1 \}$. In this section we obtain a refined invariance result in the special case when the Hamiltonian
\[H \co \R \times (P \times \R) \to \R\]
only depends on the $t$-coordinate, and when its Hamiltonian vector field $X_H=-\sigma(t)\partial_z$ moreover satisfies
\begin{itemize}
\item $\sigma(t),\sigma'(t) \ge 0$,
\item $\sigma(t)$ has support in $\{ t \ge 1 \}$, and
\item $\sigma(t)$ is constant in $\{t \ge N\}$, for some $N \ge 1$.
\end{itemize}
Observe that this Hamiltonian isotopy fixes the hypersurfaces $\{t\} \times (P \times \R)$, where it acts by some reparametrisation of the negative Reeb flow. For a generic choice of $\sigma(t)$, it is the case that $L \cap \phi^1_H(L')$ consists of transverse double-points.

\begin{remark}
Let $H$ be a Hamiltonian as above and use $\Lambda_s'$ to denote the Legendrian submanifold of which $\phi^s_H(L')$ is an exact Lagrangian filling. Given that $\phi^s_H$ is non-trivial, for $s > 0$ sufficiently large there are no Reeb chords from $\Lambda$ to $\Lambda_s'$. Moreover, there is a natural identification
\[CF_\bullet(L,L')\simeq CF_\bullet(L,\phi^s_H(L'))=CF_\bullet^0(L,\phi^s_H(L'))\]
of graded vector spaces.

The co-complex associated to $(CF_\bullet^0(L,\phi^s_H(L'),\partial_0)$ is related to the version of the wrapped Floer cohomology complex as defined in \cite{OnWrapped}, \cite{SympGeomCot}. However, one technical difference is that the latter versions are defined using moduli-spaces of solutions to a Cauchy-Riemann equation with a perturbation-term depending on a Hamiltonian vector field.
\end{remark}

Now consider the exact Lagrangian cylinder
\[V := \phi^1_H(\R \times \Lambda')\]
and, for each $s\ge 0$, the corresponding filling $L' \odot_s V$ of $\Lambda''$.

Observe that $V$ satisfies
\[\pi_P(V)=\Pi_{\mathrm{Lag}}(\Lambda'')=\Pi_{\mathrm{Lag}}(\Lambda')\]
and that the self-intersections of $(\R \times \Lambda) \cup V$ are transverse double-points corresponding to a subset of the Reeb chords on $\Lambda \cup \Lambda'$ starting on $\Lambda$ and ending on $\Lambda'$. Moreover, there is a natural identification
\[CF_\bullet(L,L') \simeq CF_\bullet(L,L' \odot_s V)\]
of graded vector spaces.

The below result shows that this identification may be assumed to hold on the level of complexes as well, given that we choose the almost complex structure with some care.

To that end, we consider a compatible almost complex structure $J_P$ on $(P,d\theta)$ and let $\widetilde{J}_P$ denote its cylindrical lift. We let $J$ be a compatible almost complex structure on $\R \times (P \times \R)$ coinciding with $\widetilde{J}_P$ in the set $\{ t \ge 1 \}$. The following proposition is also a key step in the proof of Theorem \ref{thm:twocopy} below.

\begin{proposition}
\label{prop:transfer}
Let $J_P$ be a regular compatible almost complex structure on $P$ that is integrable in a neighbourhood of the double-points of $\Pi_{\mathrm{Lag}}(\Lambda \cup \Lambda')$ and consider an almost complex structure $J$ as above, that moreover satisfies \ref{ra} in a neighbourhood of the double-points $L \cap L' \subset \{t < 0 \}$. After an arbitrarily small compactly supported perturbation of $J$, and for $s>0$ sufficiently large, we may suppose that
\begin{enumerate}
\item $J$ is regular for the moduli spaces in the definition of $(CF_\bullet(L,L'),\partial)$ and $(CF_\bullet(L,L' \odot_s V),\partial''_s)$.
\item The transfer map
\[ \Phi_V \co (CF_\bullet(L,L'),\partial) \to (CF_\bullet(L,L' \odot_s V),\partial''_s)\]
is the identity map with respect to the natural identification of the generators.
\end{enumerate}
\end{proposition}

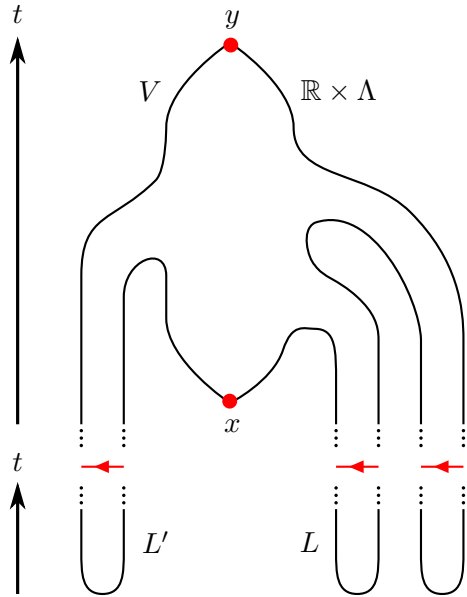
\begin{figure}[htp]
\centering
\labellist
\pinlabel $x$ at 84 64
\pinlabel $y$ at 84 217
\pinlabel $t$ at 3 219
\pinlabel $t$ at 3 50
\pinlabel $\R \times \Lambda$ at 124 190
\pinlabel $V$ at 53 190
\pinlabel $L'$ at 55 20
\pinlabel $L$ at 113 20
\endlabellist
\includegraphics{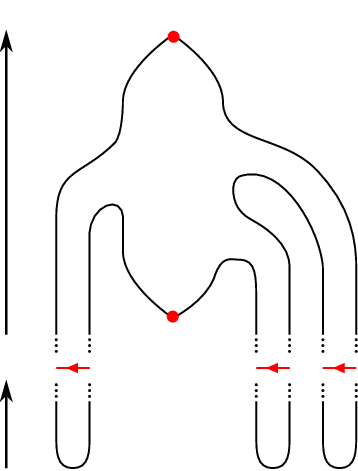}
\caption{A possible pseudo-holomorphic building in $\lim_{s \to +\infty} \mathcal{M}_{y;x}(L \cup (L' \odot_s V);J)$ consisting of a top and a bottom level.}
\label{fig:a}
\end{figure}

\begin{figure}[htp]
\centering
\labellist
\pinlabel $y$ at 85 170
\pinlabel $t$ at 5 172
\pinlabel $t$ at 5 58
\pinlabel $\R \times \Lambda$ at 124 145
\pinlabel $V$ at 56 145
\pinlabel $L'$ at 70 30
\pinlabel $L$ at 132 38
\pinlabel $x$ at 173 145
\pinlabel $L'$ at 18 30
\pinlabel $\R \times \Lambda$ at 140 110
\pinlabel $V$ at 195 110
\endlabellist
\includegraphics{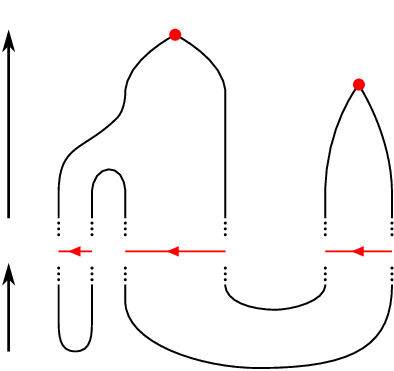}
\caption{A possible pseudo-holomorphic building in $\lim_{s \to +\infty} \mathcal{M}_{y;x}(L \cup (L' \odot_s V);J)$ consisting of a top and a bottom level.} \label{fig:b}
\end{figure}
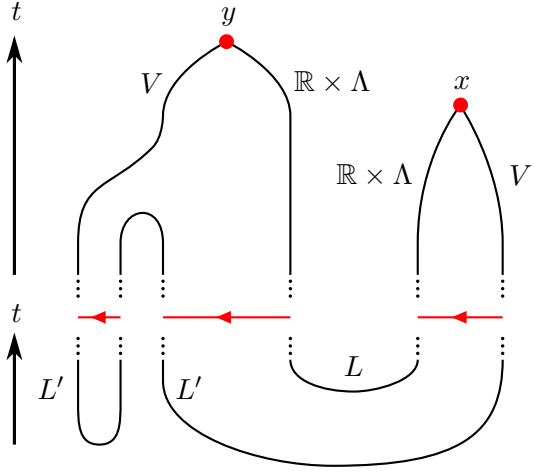

\begin{proof}
(1): A $J$-holomorphic disc having a single positive boundary-puncture asymptotic to a Reeb chord is somewhere injective in a set $\{ t \ge M \}$, as follows by its asymptotical properties. By a standard transversality argument \cite[Chapter 3]{JholCurves}, it follows that the moduli spaces consisting of such discs may be assumed to be transversely cut out after a perturbation as above.

For a regular almost complex structure $J_P$, Theorem \ref{thm:lift} implies that the $\widetilde{J}_P$-holomorphic curves with boundary on $\R \times (\Lambda \cup \Lambda')$ also can be assumed to be regular.

The moduli spaces of $J$-holomorphic strips, of which at least one boundary puncture maps to a double-point in $L \cap L'$, may be assumed to be transversely cut out by the transversality argument in \cite[Proposition 2.3]{ContHomP}. Here we have used that \ref{ra} holds in a neighbourhood of these double-points, and that the two punctures of the strip correspond to different generators, as follows by the formula for the symplectic area.

What remains is to show the transversality of the moduli spaces
\[\mathcal{M}_s:=\mathcal{M}_{y;x}(L \cup (L' \odot_s V);J)\]
of rigid $J$-holomorphic strips whose boundary-punctures both map to double-points $x, y \in (\R \times \Lambda) \cap V$.

To that end, we will consider the limit of the boundary conditions $L \cup (L' \odot_s V)$ as $s \to +\infty$. This limit amounts to stretching the neck along the contact-type hypersurface $\{t=2\}$, see \cite[Section 3.4]{CompSFT}. Furthermore, by the Gromov-Hofer compactness in \cite{CompSFT}, it follows that the solutions in $\mathcal{M}_s$ converge to so-called pseudo-holomorphic buildings consisting of the following levels.
\begin{itemize}
\item A top level consisting of $\widetilde{J}_P$-holomorphic discs with boundary on $(\R \times \Lambda) \cup V$ and boundary-punctures mapping to double-points and Reeb chords. This level is non-empty and can either consist of a single (possibly broken) disc, as shown in Figure \ref{fig:a}, or of two (possibly broken) discs, as shown in Figure \ref{fig:b}.
\item Possibly several middle levels consisting of $\widetilde{J}_P$-holomorphic discs with boundary on the cylindrical Lagrangian submanifold $\R \times (\Lambda \cup \Lambda')$.
\item A bottom level consisting of (possibly broken) $J$-holomorphic discs with boundary on $L \cup L'$. Furthermore, each disc has either exactly one positive puncture, as shown in Figure \ref{fig:a}, or exactly two, as shown in Figure \ref{fig:b}. In any case, each disc is somewhere injective.
\end{itemize}
By the index formula it can be computed that the expected dimensions of all the moduli spaces of pseudo-holomorphic discs in the building sum to $0-k$, where $k \ge 0$ is the total number of nodes of the broken discs. By a broken disc, we here mean broken in the Floer-sense, i.e.~living in the boundary of a moduli-space compactified as in \cite[Section 2.3]{ContHomR}. A broken disc thus consists of several pseudo-holomorphic discs joined at nodes, where each node corresponds to a pair of boundary-punctures that are mapped to the same double-point.

Note that, in order to see that the only possible configurations of the top layer are the two (possibly broken) configurations described above, we have used the exactness of the boundary condition. Namely, together with Stoke's theorem, it follows that there can be no (broken) disc, all whose punctures are negative punctures asymptotic to Reeb chords. 

By the above reasoning, the discs in the bottom level may be assumed to be transversely cut out. Furthermore, since $J_P$ is regular, Theorem \ref{thm:lift} implies that the solutions in the middle levels are transversely cut out. Using Lemma \ref{lem:trans2} it follows that the top levels are transversely cut out as well. Since all the moduli spaces in the building are transversely cut out, holomorphic gluing finally shows that $\mathcal{M}_s$ is transversely cut out, given that $s>0$ is sufficiently large. This finishes the claim.

As a side remark, we make the following observations. It follows that each disc in the building must live in a moduli space of non-negative expected dimension. The additivity of the dimension thus implies that all these moduli spaces are zero-dimensional, and thus no disc is broken in the Floer sense. It also follows that the middle-levels all are trivial strips over Reeb chords, and can thus be ignored.

Finally, it can be seen that there is no rigid disc as shown on the right in the top level of Figure \ref{fig:b}. To that end, observe that such a disc would project to a $J_P$-holomorphic disc having boundary on $\Pi_{\mathrm{Lag}}(\Lambda \cup \Lambda')$ of negative expected dimension, which therefore is constant. However, no disc as shown on the right in the top level of Figure \ref{fig:b} is contained in a plane of the form $\R \times (\{q\} \times \R)$.

(2): Since $\pi_P(V)=\Pi_{\mathrm{Lag}}(\Lambda')$, the $(\widetilde{J}_P,J_P)$-holomorphic projection $\pi_P$ induces maps
\begin{gather*}
\mathcal{M}_{f;\mathbf{a},c,\mathbf{b}}( (\R \times \Lambda) \cup V;\widetilde{J}_P) \to \mathcal{M}_{f;\mathbf{a},c,\mathbf{b}}( \Pi_{\mathrm{Lag}}(\Lambda \cup \Lambda');J_P) \\
\widetilde{u} \mapsto \pi_P \circ \widetilde{u}
\end{gather*}
between moduli spaces, where $c$ is a Reeb chord starting on $\Lambda$ and ending on $\Lambda'$, $f$ is either a double-point or a mixed Reeb chord, and where $\mathbf{a}=a_1\cdot\hdots\cdot a_{m_1}$ and $\mathbf{b}=b_1\cdot\hdots\cdot b_{m_1}$ are words of Reeb chords on $\Lambda'$ and $\Lambda$, respectively.

By definition $|f|-|c|-|\mathbf{a}|-|\mathbf{b}|=0$ holds for the moduli spaces contributing to $\Phi_V$. Consequently, their projections are $J_P$-holomorphic discs having exactly one positive puncture and negative expected dimension, as follows by the index computations in Lemma \ref{lem:trans2}. By the regularity of $J_P$, these latter discs must be constant, which implies that the corresponding $\widetilde{J}_P$-holomorphic discs are actually strips contained entirely in planes of the form $\R \times (\{q\} \times \R)$.

From this it follows that $\Phi_V$ has the required form. Here, part (1) of Lemma \ref{lem:trans2} has been used to infer that the above moduli spaces contributing to $\Phi_V$ are transversely cut out.

Recall that, we might have had to perturb $\widetilde{J}_P$ by a compactly supported perturbation in step (1) to achieve transversality for some of the moduli spaces under consideration, and that there is no guarantee that this perturbation is cylindrical. However, the computations in part (2) are still valid, given that the perturbation is sufficiently small. Namely, one may assume that there is a bijection between the moduli spaces in the definition of $\Phi_V$ before and after such a perturbation.
\end{proof}

\subsection{Consequences of the invariance}

An immediate consequence of the invariance theorem is the acyclicity of the wrapped Floer homology in the current setting.
\begin{proposition}
\label{prop:vanish}
Let $L,L'$ be exact Lagrangian fillings inside the symplectisation of a contactisation. It follows that $HF_\bullet(L,L')=0$.
\end{proposition}
\begin{proof}
Using the negative Reeb flow $-\partial_z$ one can isotope $L'$ to an exact Lagrangian filling $L''$ for which $CF_\bullet(L,L'')=0$. Since this is a Hamiltonian flow, the claim now follows from the invariance theorem.
\end{proof}

A similar argument shows an analogous property for the linearised Legendrian contact cohomology of Legendrian links in certain positions.

\begin{proposition}
\label{prop:vanishlch}
Let $L,L' \subset \R \times (P \times \R)$ be exact Lagrangian fillings of $\Lambda$ and $\Lambda'$, respectively, which induce an augmentation $\epsilon_{L \cup L'}$ of the Chekanov-Eliashberg algebra of the link $\Lambda \cup \Lambda'$. If all mixed Reeb-chords on $\Lambda \cup \Lambda'$ start on $\Lambda$, it follows that
\[HCL^\bullet(\Lambda, \Lambda';\epsilon_L,\epsilon_{L'})=0,\]
or equivalently,
\[HCL^\bullet(\Lambda \cup \Lambda';\epsilon_{L \cup L'})=HCL^\bullet(\Lambda;\epsilon_L) \oplus HCL^\bullet(\Lambda';\epsilon_{L'}).\]
\end{proposition}
\begin{proof}
There is a Lagrangian filling $L''$ isotopic to $L'$ by a compactly supported Hamiltonian isotopy which satisfies $CF^0_\bullet(L,L'')=0$. It immediately follows that
\[(CF_\bullet(L,L''),\partial'')=(CF^\infty_\bullet(L,L''),\partial_\infty)=(CL^{\bullet-2}(\Lambda, \Lambda'),d_{\epsilon_L,\epsilon_{L'}}).\]
Proposition \ref{prop:vanish} implies that the above complexes are acyclic.
\end{proof}

In particular, the assumptions of the previous proposition are fulfilled for a link consisting of a Legendrian submanifold $\Lambda \subset P \times \R$ admitting a filling inside the symplectisation, together with a copy of $\Lambda$ translated sufficiently far in the positive $z$-direction. This shows that $\Lambda$ satisfies the requirements for the existence of the duality long exact sequence in \cite{DualityLeg}, given that we use the augmentation induced by the filling.

\begin{example}
\label{ex:nondisp}
\begin{enumerate}
\item The fact that the augmentation is induced by a filling \emph{inside the symplectisation} is crucial for the above to hold. For instance, let $M$ be a compact manifold with boundary and let $\Lambda$ be the zero-section of $J^1(\partial M)$. The closed Legendrian submanifold $\Lambda$ has a filling $L$ consisting of the zero-section of the symplectic manifold $T^*M$, where the latter is considered as an exact symplectic manifold having a convex cylindrical end over the contact manifold $J^1(\partial M)$.

We consider a translation $\Lambda'$ of $\Lambda$ sufficiently far in the positive $z$-direction, and let $L'$ denote a filling of $\Lambda'$ that is Hamiltonian isotopic to $L$. Observe that, since neither $\Lambda$ nor $\Lambda'$ have any Reeb chords, these fillings necessarily induce the trivial augmentations.

The results in Section \ref{sec:twocopy} can be used to show that the corresponding linearised Legendrian contact cohomology satisfies
\[ HCL^\bullet(\Lambda, \Lambda';\epsilon_L,\epsilon_{L'}) = H_{\bullet+1}(\partial M;\Z_2),\]
which clearly is non-zero.
\item Let $\Lambda_1 \subset J^1(N)$ be the zero-section, and let $\Lambda_2$ be a copy of the zero-section shifted in the positive $z$-direction. It is readily seen that the Legendrian link $\Lambda:=\Lambda_1 \cup \Lambda_2$ has an exact Lagrangian filling $L \subset \R \times J^1(N)$ diffeomorphic to $\R \times N$. The Legendrian submanifold $\Lambda$ is \emph{not} horizontally displaceable, since a calculation as above shows that
\[ HCL^\bullet(\Lambda, \Lambda';\epsilon_0,\epsilon_0) = \oplus_{i=1}^4 H_{\bullet+1}(N ;\Z_2) \neq 0\]
for the trivial augmentation $\epsilon_0$, where $\Lambda'$ again is a copy of $\Lambda$ translated sufficiently far in the positive $z$-direction. In contrast to this, it follows from Proposition \ref{prop:vanishlch} that
\[ HCL^\bullet(\Lambda, \Lambda';\epsilon_L,\epsilon_{L'}) =0,\]
where $L'$ denotes a suitable translation of the filling $L$. Consequently, \cite[Theorem 1.1]{DualityLeg} can still be applied to $\Lambda$, given that we use the augmentation induced by the filling.
\end{enumerate}
\end{example}

\section{Applications of Theorem \ref{thm:lift}}
\label{sec:cor}

In this section we use Theorem \ref{thm:lift} to show that the analysis in \cite{DualityLeg} made for pseudo-holomorphic polygons in $P$ with boundary on an exact Lagrangian immersion carries over to the symplectisation. 

Let $(f,g,J_P)$ be a triple adjusted to $\Pi_{\mathrm{Lag}}(\Lambda)$ as defined in Section \ref{sec:acs}, where $f\co\Lambda \to (0,1/2]$ is a Morse function, $g$ is a Riemannian metric on $\Lambda$, and $J_P$ is a compatible almost complex structure on $P$ induced by the metric $g$. Suppose that $\Pi_{\mathrm{Lag}}(\Lambda')$ is an exact Lagrangian immersion that is $C^\infty$-close to $\Pi_{\mathrm{Lag}}(\Lambda)$. Furthermore, we assume that $\Pi_{\mathrm{Lag}}(\Lambda')$ is identified with the image of the section $df$ under the symplectic immersion of the co-disc bundle $(D^*\Lambda,d\theta_\Lambda)$ used in the construction of $J_P$ (see Section \ref{sec:acsmetric}).

In \cite{DualityLeg} it is shown that the rigid $J_P$-holomorphic polygons having one positive puncture and boundary on $\Pi_{\mathrm{Lag}}(\Lambda \cup \Lambda')$ correspond to $J_P$-holomorphic polygons having boundary on $\Pi_{\mathrm{Lag}}(\Lambda)$ together with gradient flow-lines on $\Lambda$. In Section \ref{sec:twocopy} below we recall these results.

In the case when $\widetilde{J}_P$ is the cylindrical lift of $J_P$ to the symplectisation, we can use the lifting result in Theorem \ref{thm:lift} to obtain an analogous result for the corresponding $\widetilde{J}_P$-holomorphic discs inside the symplectisation having boundary on $(\R \times \Lambda) \cup (\R \times \Lambda')$.

In order to apply Theorem \ref{thm:lift}, recall that $J_P$ needs to be integrable near the double-points of $\Pi_{\mathrm{Lag}}(\Lambda \cap \Lambda')$. For the double-points corresponding to double-points of $\Pi_{\mathrm{Lag}}(\Lambda)$ this automatically follows from the fact that $J_P$ is induced by a metric on $\Lambda$ as in Section  \ref{sec:acsmetric}. For the double-points corresponding to $\Crit(f)$, this can be achieved by choosing the metric $g$ in the construction of $J_P$ to be flat near $\Crit(f)\subset \Lambda$.

Finally, we use these results to complete the proof of Seidel's isomorphism outlined in \cite{RationalSFT2}.

\subsection{The Legendrian contact cohomology of the Legendrian two-copy link}
\label{sec:twocopy}
For an $n$-dimensional Legendrian submanifold $\Lambda \subset P \times \R$ we define the following Legendrian two-copy links. Let $(f,g,J_P)$ be a triple adjusted to $\Pi_{\mathrm{Lag}}(\Lambda)$.

The construction of the almost complex structure $J_P$ induced by the metric $g$ in Section \ref{sec:acsmetric} provides the choice of a symplectic immersion of $(D^*\Lambda,d\theta_\Lambda)$ into $P$, by which the zero-section is mapped to $\Pi_{\mathrm{Lag}}(\Lambda)$. This immersion lifts to a contact-form preserving diffeomorphism from a neighbourhood of the zero-section in $(J^1(\Lambda),dz+\theta_\Lambda)$ to a neighbourhood of $\Lambda \subset P \times \R$, which moreover maps the zero-section to $\Lambda$. We consider the following Legendrian submanifolds, where the constants $\epsilon>\eta>0$ will become useful in Section \ref{sec:floercomp} below.
\begin{itemize}
\item Let $\Lambda_+ \subset P \times \R$ be the push-off of $\Lambda$ corresponding to $(\eta^2df,\epsilon^2+ \eta^2 f) \subset J^1(\Lambda)$ under the above identification, where $\epsilon>\eta>0$ are chosen to be sufficiently small. In particular we assume that $\epsilon^2+\eta^2$ is smaller than the shortest Reeb chord on $\Lambda$.
\item Let $\Lambda_- \subset P \times \R$ be the push-off of $\Lambda$ corresponding to $(\eta^2df,-\epsilon^{1/2}+ \eta^2 f) \subset J^1(\Lambda)$ under the above identification, where $\epsilon>\eta>0$ are chosen to be sufficiently small. In particular, we assume that $\epsilon^{1/2}-\eta^2>0$ is smaller than the shortest Reeb chord on $\Lambda$.
\item Let $\Lambda_\infty \subset P \times \R$ be a copy of $\Lambda_+$ translated sufficiently far in the positive $z$-direction, so that all mixed Reeb chords on $\Lambda \cup \Lambda_\infty$ start on $\Lambda$. 
\end{itemize}

\subsubsection{Generalised pseudo-holomorphic discs} We recall that a \emph{generalised pseudo-holomorphic disc} in $P$ consists of the following data. Let $c$ be a critical point of $f$ and let
\[ u \co (D^2,\partial D^2) \to (P,\Pi_{\mathrm{Lag}}(\Lambda)) \]
be a pseudo-holomorphic polygon having boundary on $\Pi_{\mathrm{Lag}}(\Lambda)$ and boundary punctures mapping to double-points, together with an additional marked point $p_f \in \partial D^2$. We require $u(p_f)$ to be connected to $c$ via a flow-line of $\nabla f$ (we allow $c=u(p_f)$). In the case when $c$ is connected to $u(p_f)$ by the positive (respectively negative) gradient flow, we say that $p_f$ is a positive (respectively negative) Morse puncture.

We refer to \cite{DualityLeg} for the expected dimension and transversality results for these moduli spaces.

\subsubsection{Computations of the complexes for the different two-copies} We assume that the Chekanov-Eliashberg algebra of $\Lambda$ has an augmentation $\epsilon$. We may assume that $\Lambda_\pm$ are sufficiently $C^\infty$-close to $\Lambda$ and hence that the Chekanov-Eliashberg algebra of $\Lambda_i$ for $i=+,-,\infty$ coincide with the Chekanov-Eliashberg algebra of $\Lambda$ (see e.g.~\cite[Lemma 4.13]{OrientLeg}). In particular, $\epsilon$ induces an augmentation of each component of the link $\Lambda \cup \Lambda_i$, and thus of the link itself.

We use $(CL^\bullet(\Lambda,\Lambda_i),d_i)$ to denote the linearised Legendrian contact cohomology complex generated by Reeb chords starting on $\Lambda$ and ending on $\Lambda_i$, where the complex is linearised using the augmentation $\epsilon$.

Choosing canonically defined capping paths, it can be shown that
\begin{eqnarray*}
& & CL^\bullet(\Lambda,\Lambda_-)=CL^\bullet(\Lambda), \\
& & CL^\bullet(\Lambda,\Lambda_+)= C_{\mathrm{Morse}}^{\bullet+1}(f) \oplus CL^\bullet(\Lambda), \\
& & CL^\bullet(\Lambda,\Lambda_\infty)=CL_{n-2-\bullet}(\Lambda) \oplus C_{\mathrm{Morse}}^{\bullet+1}(f) \oplus CL^\bullet(\Lambda),
\end{eqnarray*}
where we recall that $n=\dim \Lambda$. To that end, observe that a mixed Reeb chord either corresponds to a Reeb chord on $\Lambda$, or to a critical point of the Morse function $f$. We refer to \cite[Section 3.1]{DualityLeg} for more details.

Since the co-differential is action-increasing, by comparing the length of the Reeb chords in the different summands above, one concludes that the differentials are of the following form with respect to the above decompositions.
\begin{eqnarray*}
& & d_- = d_q, \\
& & d_+ = \begin{pmatrix}
d_f & 0\\
\sigma & d_q
\end{pmatrix}, \\
& & d_\infty = \begin{pmatrix}
d_p & 0 & 0 \\
\rho & d_f & 0 \\
\eta & \sigma & d_q
\end{pmatrix}.
\end{eqnarray*}

\begin{remark} Note that the exact Lagrangian immersion $\Pi_{\mathrm{Lag}}(\Lambda \cup \Lambda_i)$ is the same for $i=-,+,\infty$ and, hence, the set of pseudo-holomorphic polygons in $P$ with boundary on $\Pi_{\mathrm{Lag}}(\Lambda \cup \Lambda_i)$ is independent of $i$. However, the notion of being a positive or a negative puncture does depend on $i$.
\end{remark}

By \cite[Lemma 6.5]{DualityLeg}, we may assume that $(f,g)$ is a Morse-Smale pair and that $J_P$ is regular. Furthermore, according to \cite[Theorem 3.6]{DualityLeg}, for generic such triples $(f,g,J_P)$ adjusted to $\Pi_{\mathrm{Lag}}(\Lambda)$, and for $\eta>0$ small enough, the above co-differentials can be defined by counting the following (generalised) $J_P$-holomorphic discs. 
\begin{itemize}
\label{list:multiple}
\item $d_p=\partial_\epsilon$ is the differential of the linearised Legendrian contact homology complex of $\Lambda$ with respect to the augmentation $\epsilon$.
\item $d_q=d_\epsilon$ is the differential on the linearised Legendrian contact cohomology complex with respect to the augmentation $\epsilon$.
\item The differential $d_f$ is the differential on the Morse co-complex (i.e. counting positive gradient flow-lines).
\item The map $\sigma$ counts rigid generalised pseudo-holomorphic discs with boundary on $\Pi_{\mathrm{Lag}}(\Lambda)$, one positive puncture, and one negative Morse-puncture at a critical point of $f$.
\item The map $\eta$ counts rigid pseudo-holomorphic discs with boundary on $\Pi_{\mathrm{Lag}}(\Lambda)$ and two positive punctures (after a domain-dependent perturbation of the boundary condition).
\item The map $\rho$ counts generalised pseudo-holomorphic discs with boundary on $\Pi_{\mathrm{Lag}}(\Lambda)$ and two positive punctures, of which one is a Morse-puncture at a critical point of $f$.
\end{itemize}

Note that requiring $g$ to be flat in neighbourhoods of the critical points of $f$ does not impose any restriction here. In the case when $\widetilde{J}_P$ is the cylindrical lift of $J_P$, Theorem \ref{thm:lift} thus applies to give that the above description of the complexes holds for the version of Legendrian contact homology defined in terms of $\widetilde{J}_P$-holomorphic discs in the symplectisation as well. To that end, for such a choice of $g$, the induced almost complex structure $J_P$ becomes integrable in a neighbourhood of the double-points of $\Pi_{\mathrm{Lag}}(\Lambda \cup \Lambda_i)$. 

\subsection{The wrapped Floer homology of the two-copy of a filling}
\label{sec:floercomp}

Let $L \subset \R \times (P \times \R)$ be an $(n+1)$-dimensional exact Lagrangian filling of the Legendrian submanifold $\Lambda$. We assume that $L$ is cylindrical in the set $\{ t \ge -2 \}$.

\subsubsection{Construction of the Morse functions $F_\pm$ on $L$}
\label{sec:morse}
Start by fixing a positive Morse function
\[f \co \Lambda \to (0,1/2],\]
and a smooth cut-off function $\rho \co \R_{>0} \to [0,1]$ satisfying $\rho'(s) \ge 0$, $\rho|_{(0,e^{-1}]}=0$, and $\rho|_{[e^{-1/2},+\infty)}=1$. For each $0<\eta <1$ we define the smooth function
\begin{gather*}
f_\eta \co \R_{> 0} \times \Lambda \to \R_{>0}, \\
(s,p) \mapsto f_\eta(s,p):= \eta^2 \rho(e^{-1/2}s) f(p),
\end{gather*}
which thus coincides with $\eta^2f$ on $\{s \ge e^0\}$ and vanishes on $\{s \le e^{-1/2}\}$. Also, for each $1>\epsilon>\eta>0$ and $B>0$, we define the smooth cut-off functions
\begin{eqnarray*}
\sigma_\epsilon(s) & := &  \epsilon^5+(\epsilon^2-\epsilon^5)\rho(s),\\
\widetilde{\sigma}_\epsilon(s) & := &  \epsilon^5+(\epsilon^2-\epsilon^5)\rho(s)+(\epsilon^{1/2}-\epsilon^2)\rho(e^{-(B+3/2)}s).
\end{eqnarray*}
It follows that $\sigma_\epsilon|_{(0,e^{-1}]}=\epsilon^5$, $\sigma_\epsilon|_{[e^{-1/2},+\infty)}=\epsilon^2$, while $\widetilde{\sigma}_\epsilon|_{(0,e^{B+1/2}]}=\sigma_\epsilon$, and $\widetilde{\sigma}_\epsilon|_{[e^{B+1},+\infty)}=\epsilon^{1/2}$.

Consider a fixed Morse function $G$ on $L \cap \{ t \le -1 \}$ given by $e^t|_L$ on
\[L \cap \{-2 \le t \le -1 \}=[-2,-1] \times \Lambda \subset \R \times (P \times \R).\]
We define
\[F^{\eta,\epsilon}_+ \co L \to \R\]
to be the Morse function coinciding with $\epsilon^5 G$ on $L \cap \{ t \le -1 \}$ and with $\sigma_\epsilon(s) s+sf_\eta$ on
\[L \cap \{t  \ge -1\}=[-1,+\infty) \times \Lambda \subset \R \times (P \times \R),\]
where we have set $s:=e^t|_{L \cap \{t  \ge -1\}}$.

Let $0 < A < B$ be fixed. We also consider the Morse function
\[F^{\eta,\epsilon}_- \co L \to \R\]
coinciding with $F^{\eta,\epsilon}_+$ on
\[\overline{L}:=L \cap \{ t \le A-1\}\]
and which is given by $\widetilde{\sigma}_\epsilon(s) \alpha(s)+sf_\eta$ on
\[L \cap \{t \ge A-1\}=[A-1,+\infty) \times \Lambda \subset \R \times (P \times \R).\]
Here $0< A < B$ have been chosen sufficiently large, and $\alpha \co \R_{>0} \to \R$ is a smooth function satisfying
\begin{itemize}
\item $\alpha(s)=s$ for $s \le e^{A-1}$, $\alpha(s)>0$ for $e^A \le s \le e^B$, and $\alpha(e^{B+1/2})=0$,
\item $\alpha'(e^A)=0$, $\alpha'(e^B)<-1/2$, and $\alpha'(s)=-1$ for $s \ge e^{B+1}$, and
\item $\alpha''(s) \le 0$ for all $s$, and $\alpha''(s)<0$ for $e^A \le s \le e^B$.
\end{itemize}
We will sometimes use $F_\pm$ to denote $F^{\eta,\epsilon}_\pm$.

The critical points of $F_+$ correspond to the critical points of $G$ and are all contained in $\overline{L}$.

The critical points of $F_-$ consist of critical points inside $\overline{L}$, which hence correspond to the critical points of $F_+$, together with critical points inside the set $L \cap \{A < t < B \}$. The latter critical points are determined by the equation
\[ (\epsilon^2\alpha'(s)+ \eta^2f(p))ds+s \eta^2 df(p)=0.\]
Since $\alpha''(s)<0$ holds in this set, it follows that these critical points are non-degenerate and correspond bijectively to critical points of $f$. Moreover, the Morse index of such a critical point is one greater than the Morse index of the corresponding critical point of $f$.

We let $g$ be a choice of a metric on $\Lambda$ for which $(f,g)$ is Morse-Smale and let $(C^{\bullet-1}_{\mathrm{Morse}}(f),d_f)$ be the induced Morse co-complex, whose differential thus counts \emph{positive} gradient flow lines.

Also, consider a Riemannian metric on $L$ coinciding with a perturbation of the product metric $dt \otimes dt +g$ in the set $\{t \ge 0\}$, such that this metric together with a small perturbation of $F_-$ is Morse-Smale. According to Lemma \ref{lem:morse} we may assume that the corresponding Morse co-complex
\[ (C^\bullet_{\mathrm{Morse}}(F_-),d_{F_-})=(C^\bullet_{\mathrm{Morse}}(F_+) \oplus C^{\bullet-1}_{\mathrm{Morse}}(f),d_{F_-}),\]
has differential given by
\[d_{F_-}= \begin{pmatrix} d_{F_+} & 0 \\
\Gamma & d_f
\end{pmatrix}.\]

Finally, the Morse cohomology groups associated to $F_\pm$ are given by
\begin{gather*}
H^\bullet_{\mathrm{Morse}}(F_+;\Z_2)=H^\bullet(L;\Z_2)=H_{(n+1)-\bullet}(\overline{L},\partial \overline{L};\Z_2),\\
H^\bullet_{\mathrm{Morse}}(F_-;\Z_2)=H^\bullet(\overline{L},\partial \overline{L};\Z_2)=H_{(n+1)-\bullet}(L;\Z_2).
\end{gather*}
The equalities on the left follow from the fact the differential in Morse cohomology counts the positive gradient flow-lines, while the negative gradients of $F_+$ and $F_-$ points inwards and outwards of $L \cap \{t \le B \}$, respectively. The equalities on the right follow by the Poincar\'{e} duality for manifolds with boundary.

\subsubsection{Construction of the push-offs $L_\pm$ of $L$} We will use the above Morse functions $F_\pm$ to construct certain exact Lagrangian fillings corresponding to push-offs of $L$.

As above, we fix a symplectic immersion of $(D^*\Lambda,d\theta_\Lambda)$ into $P$ extending the immersion $\iota$, as constructed in Section \ref{sec:acsmetric}. This symplectic immersion lifts to a contact-form preserving diffeomorphism
\[\varphi \co (U,dz+\theta_\Lambda) \to (V,\lambda),\:\:\lambda=dz+\theta, \]
from a neighbourhood $U \subset J^1(\Lambda)$ of the zero-section to a neighbourhood $V \subset P \times \R$ of $\Lambda$. Furthermore, $\varphi$ maps the zero-section to $\Lambda$.

It immediately follows that the map
\[(\id,\varphi^{-1}) \co ([-2,+\infty) \times V,d(e^t\lambda)) \to ([-2,+\infty) \times U,d(e^t(dz+\theta_\Lambda)))\]
is an exact symplectomorphism identifying a neighbourhood of $L \cap \{t \ge -2\}$ with a neighbourhood of the cylinder over the zero-section in the symplectisation of $(J^1(\Lambda),dz+\theta_\Lambda)$.

Furthermore, the symplectisation of $(U,dz+\theta_\Lambda)$ is symplectomorphic to a neighbourhood of the zero-section of the cotangent bundle of $\R_{>0} \times \Lambda$ via the (non-exact) symplectomorphism
\begin{gather*}
\psi \co (\R \times J^1(\Lambda),d(e^t(dz+\theta_\Lambda))) \to (T^*(\R_{>0} \times \Lambda),d\theta_{\R_{>0} \times \Lambda}),\\
(t,(\mathbf{q},\mathbf{p},z)) \mapsto ((e^t,\mathbf{q}),(z,e^t\mathbf{p})).
\end{gather*}

Using the Weinstein Lagrangian neighbourhood theorem together with the above symplectomorphisms, we can construct a symplectic identification of a neighbourhood of the zero-section in $(T^*L,d\theta_L)$ with a neighbourhood of $L$, such that the identification coincides with $(\psi \circ (\id,\varphi^{-1}))^{-1}$ on a neighbourhood of the zero-section of
\[(T^*([e^{-2},+\infty)\times \Lambda),d\theta_{[e^{-2},+\infty)\times \Lambda}) \subset (T^*L,\theta_L).\]

For $0<\eta<\epsilon<1$ sufficiently small, we use the above identification to construct the exact Lagrangian fillings $L^{\eta,\epsilon}_+$ and $L^{\eta,\epsilon}_-$ inside $\R \times (P \times \R)$ by letting them correspond to the graphs of $dF^{\eta,\epsilon}_+$ and $dF^{\eta,\epsilon}_-$ inside of $T^*L$, respectively.

By construction, we have
\begin{eqnarray*}
L^{\eta,\epsilon}_+ \cap \{ t \ge -1/2 \} & = & [-1/2,+\infty) \times \Lambda_+,\\
L^{\eta,\epsilon}_- \cap \{ t \ge B+1\} & = & [B+1,+\infty) \times \Lambda_-,
\end{eqnarray*}
where the Legendrian submanifold $\Lambda_+$ corresponds to the 1-jet
\[(d(\epsilon^2+ \eta^2 f),\epsilon^2+ \eta^2 f) \subset (T^*\Lambda \times \R=J^1(\Lambda),dz+\theta_\Lambda),\]
and $\Lambda_-$ corresponds to the 1-jet
\[(d(-\epsilon^{1/2}+ \eta^2 f),-\epsilon^{1/2}+ \eta^2 f) \subset (T^*\Lambda \times \R=J^1(\Lambda),dz+\theta_\Lambda),\]
under the above identification. In other words, the Legendrian ends correspond to the push-offs of $\Lambda$ as constructed in Section \ref{sec:twocopy} above.

Finally, note that the double-points of $L \cup L_{\pm}$ are in bijective correspondence with the critical points of $F_\pm$.

\subsubsection{Computation of $(CF(L,L_+),\partial)$} We are now ready to state and prove the main result of this section.
\begin{theorem}
\label{thm:twocopy}
The wrapped Floer homology complex
\[(CF_\bullet(L,L^{\eta,\epsilon}_+),\partial)=(CF_\bullet^0(L,L_+) \oplus CF_\bullet^\infty(L,L_+),\partial),\]
is given by
\begin{eqnarray*}
& & CF_\bullet^0(L,L_+)=C^{\bullet}_{\mathrm{Morse}}(F_+),\\
& & CF_\bullet^\infty(L,L_+)=C_{\mathrm{Morse}}^{\bullet-1}(f) \oplus CL^{\bullet-2}(\Lambda),
\end{eqnarray*}
where the differential is of the form
\[\partial = \begin{pmatrix}
d_{F_+} & 0 &0 \\
\gamma & d_f & 0  \\
g & \sigma & d_q
\end{pmatrix}.\]
For $0<\eta<\epsilon<1$ sufficiently small, a suitable Morse-function $f$ on $\Lambda$, and a suitable compatible almost complex structure on $\R \times (P \times \R)$, it follows that
\begin{enumerate}
\renewcommand{\labelenumi}{(\roman{enumi})}
\item $(CL^{\bullet-2}(\Lambda),d_q)$ is the linearised Legendrian contact cohomology complex for $\Lambda$ induced by the filling $L$, $\sigma$ is the map described in Section \ref{sec:twocopy}, and $d_f$ is the Morse co-differential induced by $f$.
\item $d_{F_+}$ is the Morse co-differential induced by $F_+$.
\item $\gamma$ is homotopic to $\Gamma$ defined above.
\end{enumerate}
In particular, $\Cone(\gamma)$ is isomorphic to $(C^\bullet_{\mathrm{Morse}}(F_-),d_{F_-})$.
\end{theorem}

\begin{remark}
\label{rem:conj}
The proof of Corollary \ref{cor:dual} given in \cite{RationalSFT2} depends on the conjectural analytical results (1)-(5) stated in \cite[Conjectural Lemma 4.11]{RationalSFT}. More precisely, the following statements are needed. First, the proof uses statements (1)-(3) together with a translation of the results in \cite{DualityLeg} to the version of Legendrian contact homology defined via the symplectisation. These results follow from Theorem \ref{thm:lift}, as shown in Section \ref{sec:twocopy} above. Second, the proof depends on the homotopy equivalence established in part (iii) of Theorem \ref{thm:twocopy}, which can be seen as a weaker version of (5) in the conjectural lemma.

Observe that part (5) of the conjectural lemma claims that there exists an almost complex structure for which the pseudo-holomorphic discs in the definition of $\gamma$ are in bijection with the flow-lines in the definition of $\Gamma$. Here, we only prove that the corresponding chain-maps defined by these disc-counts are homotopy equivalent.
\end{remark}

\begin{proof}
By construction $L_+$ is cylindrical in the set $\{ t \ge -1/2 \}$, $L_-$ is cylindrical in the set $\{t \ge B+1\}$, and $L_-$ coincides with $L_+$ in the set $\{ t \le A-1 \}$. It also follows that $L_-=L_+ \odot V$, where $V$ is an exact Lagrangian cylinder. Here, we may take $V:=\phi^1_H(\R \times \Lambda_+)$ for a Hamiltonian $H \co \R \times (P \times \R) \to \R$ that only depends on the $t$-coordinate, and which satisfies the assumptions in Section \ref{sec:Reebflow}.

Let $J_P$ be a regular compatible almost complex structure on $P$ that satisfies the assumptions of Theorem \ref{thm:lift}. We use $\widetilde{J}_P$ to denote its cylindrical lift. We fix a compatible almost complex structure $J$ on $\R \times (P \times \R)$ that is induced by a Riemannian metric on $L$ in some neighbourhood of $L \cap \{ t \le 0 \}$ and which coincides with $\widetilde{J}_P$ on $\{ t \ge 1 \}$. 

After an arbitrarily small compactly supported perturbation in the set $\{ t \ge B+3\}$, we assume that $J$ is regular for the moduli spaces in the definition of $(CF_\bullet(L,L^{\eta,\epsilon}_+),\partial)$.

(i): We will choose a triple $(f,g,J_P)$ adjusted to $\Pi_{\mathrm{Lag}}(\Lambda)$, where $J_P$ is regular and $(f,g)$ is Morse-Smale. As before, we will choose the metric $g$ so that $J_P$ becomes integrable in some neighbourhood of $\Pi_{\mathrm{Lag}}(\Crit(f)) \subset \Pi_{\mathrm{Lag}}(\Lambda)$.

Recall that, by construction, there are Darboux coordinates $\mathbf{x}+i\mathbf{y} \in \C^n$ near each double-point of $\Pi_{\mathrm{Lag}}(\Lambda)$ in which $J_P=i$ and the two branches of $\Pi_{\mathrm{Lag}}(\Lambda)$ correspond to the real and imaginary part, respectively. 

The results in Section \ref{sec:twocopy} may thus be assumed to hold for the linearised Legendrian cohomology of $\Lambda \cup \Lambda_\pm$ defined in terms of $\widetilde{J}_P$ and the augmentation $\epsilon_L$ of the Chekanov-Eliashberg algebra of $\Lambda$. Recall that, for $\epsilon>0$ sufficiently small, we may assume that the Chekanov-Eliashberg algebras of $\Lambda$ and $\Lambda_\pm$ coincide.

We claim that $\partial_\infty=d_{\epsilon_L,\epsilon_{L_+}}$ coincides with $d_q=d_{\epsilon_L}$ in Section \ref{sec:twocopy}. To see this, we must show that the augmentation $\epsilon_{L_+}$ induced by $L_+$ can be taken to coincide with $\epsilon_L$ induced by $L$.

To that end, we will choose the Morse function $f \co \Lambda \to (0,1/2]$ with some care near the end-points of the Reeb chords on $\Lambda$. More precisely, in the coordinates $\varphi_{q,i} \co U_{q,i} \to D^n$ on the neighbourhoods $U_{q,i} \subset \Lambda$ of the end-points of the Reeb chords, as used in the construction of the almost complex structure in Section \ref{sec:acsmetric}, we require $f$ to be affine. It follows that there are holomorphic coordinates $\mathbf{z} \in \C^n$ near each double-point of $\Pi_{\mathrm{Lag}}(\Lambda)$ in which this Lagrangian immersion coincides with $\re(\C^n) \cup \im(\C^n)$, and where $\Pi_{\mathrm{Lag}}(\Lambda_+)$ coincides with a translation $\re(\C^n) \cup \im(\C^n)+\mathbf{z}_0$.

For such a Morse function $f$, there thus exist neighbourhoods $U, U_+ \subset P$ of the double-points of $\Pi_{\mathrm{Lag}}(\Lambda)$ and $\Pi_{\mathrm{Lag}}(\Lambda_+)$, respectively, for which there is a contact-form preserving diffeomorphism
\begin{gather*}\Psi \co (U \times \R,\Lambda \cap (U \times \R) ) \to (U_+ \times \R,\Lambda_+ \cap (U_+ \times \R)),\\
(p,z) \mapsto (\psi(p),z+h(p,z)),
\end{gather*}
and where $\psi$ is the symplectomorphism given by the translation $\mathbf{z} \mapsto \mathbf{z}+\mathbf{z}_0$ in the above holomorphic Darboux coordinates. Since $\psi$ moreover is a $J_P$-biholomorphism, it follows that $(\id_\R,\Psi)$ is a $\widetilde{J}_P$-biholomorphism mapping $\R \times (\Lambda \cap (U \times \R))$ to $\R \times (\Lambda_+ \cap (U_+ \times \R))$.

After shrinking $\epsilon>\eta>0$, the two augmentations $\epsilon_L$ and $\epsilon_{L_\pm}$ of the Chekanov-Eliashberg algebra of $\Lambda_\pm$ may thus be assumed to be obtained by counts of pseudo-holomorphic discs for two almost complex structures and two Lagrangian boundary conditions that differ by an arbitrarily small compactly supported perturbation. Here we have used the existence of the above $\widetilde{J}_P$-biholomorphism together with the asymptotical properties of the solutions to infer that, outside of a compact set in the domain, two such discs are solutions to the same boundary-value problem. In particular, these two augmentations may be supposed to be equal for $\epsilon>\eta>0$ sufficiently small and a suitable Morse function $f$.

We let $(CF_\bullet(L,L_+),\partial)$ be the complex determined by $J$ and $f$ as above. The results in Section \ref{sec:twocopy} can finally be applied to give (i).

(ii): Since $J$ is induced by a metric on $L$ in a neighbourhood of $L \cap \{t \le 0\}$, part (1) of Lemma \ref{lem:mono} together with \cite[Lemma 6.11]{ContHomP} shows that, for a suitable metric on $L$ and function $F_+$, we have $\partial_0=d_{F_+}$ for $\epsilon>0$ sufficiently small.

We have concluded that the differential is of the form
\[\partial = \begin{pmatrix}
d_{F_+} & 0 &0 \\
\gamma & d_f & 0  \\
* & \sigma & d_q
\end{pmatrix}.\]

(iii): We must show that $\gamma \simeq \Gamma$, where the latter map is a component of the differential of the Morse co-complex
\begin{gather*}
(C^\bullet_{\mathrm{Morse}}(F_-),d_{F_-})=(C^\bullet_{\mathrm{Morse}}(F_+) \oplus C^{\bullet-1}_{\mathrm{Morse}}(f),d_{F_-}),\\
d_{F_-}= \begin{pmatrix} d_{F_+} & 0 \\
\Gamma & d_f
\end{pmatrix},
\end{gather*}
as constructed in Lemma \ref{lem:morse}.

Recall that $J$ coincides with a compactly supported perturbation of the cylindrical lift $\widetilde{J}_P$ in the set $\{ t \ge 1 \}$. Without loss of generality, the conditions in Proposition \ref{prop:transfer} may thus be assumed to hold, after choosing $A>0$ sufficiently big. In other words, $J$ is also regular for the moduli spaces in the definition of the complex $(CF(L,L_-),\partial_-)$ and, moreover, the obvious identification of generators induces an isomorphism of complexes. In particular, the identity $\partial_-=\partial$ holds under the canonical identification of generators.

We now consider a tame almost complex structure $J'$ satisfying the following properties.
\begin{itemize}
\item $J'$ coincides with $J$ in the set $\{ t \le 0 \} \cup \{ t \ge B+3\}$,
\item In a neighbourhood of $L \cap \{ t \le B+2 \}$, $J'$ is the push-forward of an almost complex structure on $T^*L$ induced by a metric on $L$ (see Section \ref{sec:acsmetric}).
\item The canonical projection $\pi_P$ is $(J',J_P)$-holomorphic in the set $\{ t \ge A-1\}$.
\end{itemize}

To see the existence of such an almost complex structure we argue as follows. Let $\widetilde{\iota}$ be the symplectic immersion of $(D^*\Lambda,d\theta_\Lambda)$ into $P$ that is used in the construction of the almost complex structure $J_P$ induced by the metric $g$ on $\Lambda$, as described in Section \ref{sec:acsmetric}. Recall that $\widetilde{\iota}$ restricted to the zero-section is the immersion $\iota$ of $\Pi_{\mathrm{Lag}}(\Lambda) \subset P$. Consider a (non-symplectic!) embedding
\begin{gather*}
T^*[A-2,B+3] \times T^*\Lambda \to [A-2,B+3] \times (P \times \R)\\
((t,z),(\mathbf{q},\mathbf{p})) \mapsto (t,(\widetilde{\iota}(\mathbf{q},\mathbf{p}),h(\mathbf{q},\mathbf{p})+z)),
\end{gather*}
defined in a neighbourhood of the zero-section of the domain, where $h(\mathbf{q},0)$ is given by the $z$-coordinate of the lift of $\iota$ to $\Lambda \subset P \times \R$, and where the above map moreover takes the tangent-planes $0 \oplus T(T^*\Lambda)$ along the zero-section
\[[A-2,B+3]\times\Lambda \subset T^*([A-2,B+3] \times \Lambda)\]
to the contact-planes. Observe that the zero-section is mapped to $[A-2,B+3] \times \Lambda \subset \R \times (P \times \R)$.

The construction in Section \ref{sec:acs} applied to the product metric $dt \otimes dt + g$ on $[A-2,B+3] \times \Lambda$ induces a compatible almost complex structure $J_{dt^2+g}$ on $T^*([A-2,B+3] \times \Lambda)$. Using the above (non-symplectic!) identification, this almost complex structure can be pushed forward to a neighbourhood $U$ of $[A-2,B+3] \times \Lambda \subset \R \times (P \times \R)$.

By construction, this almost complex structure $J_{dt^2+g}$ in $U$ is invariant under translations of the $t$ and $z$-coordinate, satisfies $J_{dt^2+g} \partial_t=\partial_z$, and has the property that
\[ \pi_P \co U \to \pi_P(U)\]
is $(J_{dt^2+g},J_P)$-holomorphic. Furthermore, it can be checked that $J_{dt^2+g}$ is tamed by the symplectic form on the symplectisation in some neighbourhood of $[A-2,B+3] \times \Lambda$. To construct the required tame almost complex structure $J'$, it now suffices to perform an appropriate interpolation of tame almost complex structures.

\cite[Lemma 6.11]{ContHomP} together with part (2) of Lemma \ref{lem:mono} implies that, for $\epsilon>0$ sufficiently small, we may assume that there is an identity
\[\partial'_0=d_{F_-}.\]
Here we might have to perturb the pair $(F_-,dt^2+g)$ to make it Morse-Smale. Also, we will have to use an almost complex structure $J''$ induced by this perturbed metric in a neighbourhood of $L$ in order to define $\partial'$. However, the resulting $J''$ may be assumed to be an arbitrarily small compactly supported perturbation of $J'$.

Using Lemma \ref{lem:morse} we infer that $\partial'$ is of the form
\[ \partial' = \begin{pmatrix}
d_{F_+} & 0 & 0 \\
\Gamma & d_f & 0 \\
* & * & *
\end{pmatrix}.
\]

We now consider a path $\{J_s\}_{s \in [0,1]}$ of tame almost complex structures coinciding with $J$ in the set $\{ t \le 0 \} \cup \{t \ge B+3\}$, where $J_0=J$, $J_1=J'$, and such that $\pi_P$ is $(J_s,J_P)$-holomorphic in the set $\{ t \ge A-1 \}$ for all $s$.

Using part (3) of Lemma \ref{lem:mono}, for $\epsilon>\eta>0$ sufficiently small, we may assume that every $J_s$-holomorphic strip with both punctures at generators corresponding to $C_{\mathrm{Morse}}(f)$ are contained inside $\{A-1 \le t \le B\}$. Since $\pi_P$ is $(J_s,J_P)$-holomorphic when restricted to this set, it follows that any such strip projects to a $J_P$-holomorphic strip in $P$ having boundary on $\Pi_{\mathrm{Lag}}(\Lambda \cup \Lambda')$ and being of the same index. Since $J_P$ is regular, it now follows that there can be no such $J_s$-holomorphic strips of negative index.

Furthermore, part (1) of Lemma \ref{lem:mono} implies that any $J_s$-holomorphic strip having punctures at double-points corresponding to $C_{\mathrm{Morse}}(F_+)$ are contained inside $\{ t \le 0 \}$. Since $J_s=J$ in this set, the regularity of $J$ implies that there are no such strips of negative index either.

The non-existence of the pseudo-holomorphic strips of negative index having both punctures at double-points corresponding to generators of either $C_{\mathrm{Morse}}(f)$ or $C_{\mathrm{Morse}}(F_+)$ still holds if we consider a sufficiently small perturbation $J'_s$ of the path $J_s$. Without loss of generality, we may thus assume this property to hold for a path starting at $J'_0=J$ and ending at $J'_1=J''$.

The above form of the possible $J'_s$-holomorphic discs of index $-1$ has the following consequence: The isomorphism of chain-complexes given by Proposition \ref{prop:path} applied to the path $J'_s$ is of the form
\begin{gather*}
\Phi \co (CF(L,L_-),\partial_-) \to (CF(L,L_-),\partial'),\\
\Phi=\begin{pmatrix}
\id_{C(F_+)} & 0 & 0 \\
\phi & \id_{C(f)} & 0 \\
* & * & *
\end{pmatrix}.
\end{gather*}

The chain-map property
\[\Phi\circ \partial_-=\partial' \circ \Phi,\]
can now be written as
\[\begin{pmatrix}
d_{F_+} & 0 & 0 \\
\phi \circ d_{F_+} + \gamma & d_f & 0 \\
* & * & *
\end{pmatrix}=\begin{pmatrix}
d_{F_+} & 0 & 0 \\
\Gamma+d_f \circ \phi & d_f & 0 \\
* & * & *
\end{pmatrix},
\]
which translates into the fact that there is a chain-homotopy $\Gamma \simeq \gamma$.
\end{proof}

\begin{lemma}
\label{lem:mono}
Let $L$ and $L^{\eta,\epsilon}_\pm$ be the exact Lagrangian fillings as constructed in Section \ref{sec:floercomp} and let $J_s$, $s \in [0,1]$, be a family of tame almost complex structures on $\R \times (P \times \R)$. Let $U_0$, $U_{B+2}$, and $U_{[A-1,B]}$ denote fixed compact neighbourhoods of $L \cap \{ t \le 0 \}$, $L \cap \{ t \le B+2\}$, and $L \cap \{ A-1 \le t \le B \}$, respectively. Choosing $0<\eta<\epsilon<1$ sufficiently small we may suppose that
\begin{enumerate}
\item Any $J_s$-holomorphic strip having boundary on $L \cup L^{\eta,\epsilon}_\pm$, whose both punctures are mapped to double-points corresponding to $\Crit(F^{\eta,\epsilon}_+)$, is contained inside $U_0$.
\item Any $J_s$-holomorphic strip having boundary on $L \cup L^{\eta,\epsilon}_-$, whose both punctures are mapped to double-points corresponding to $\Crit(F^{\eta,\epsilon}_-)$, is contained inside $U_{B+2}$.
\item Any $J_s$-holomorphic strip having boundary on $L \cup L^{\eta,\epsilon}_-$, whose both punctures are mapped to double-points corresponding to $\Crit(f)$, is contained inside $U_{[A-1,B]}$.
\end{enumerate}
\end{lemma}

\begin{proof}
Consider the symplectic form $\omega:=d(e^t(dz+\theta))$ on $\R \times (P \times \R)$. 

The $\omega$-area of a strip as above whose positive and negative puncture correspond to the critical points $p$ and $q$ of $F_\pm$, respectively,  is given by $F_\pm(p)-F_\pm(q)$. It immediately follows that there is a constant $E>0$ independent of $\epsilon>\eta>0$ such that any $J_s$-holomorphic strip in (1) and (2) have $\omega$-area bounded from above by $E\epsilon^5$ and $E\epsilon^2$, respectively. Similarly, it readily follows that the $J_s$-holomorphic strips in (3) may be assumed to have $\omega$-area bounded from above by $E\eta^2$, given that $\eta>0$ is sufficiently small (in particular, see the construction of $\alpha(s)$ in Section \ref{sec:morse}).

We now outline a consequence of the monotonicity property for the $\omega$-area of $J_s$-holomorphic curves without boundary \cite[Proposition 4.3.1(ii)]{SomeProp} and with boundary \cite[Proposition 4.7.2(ii)]{SomeProp}. Let $K \subset \R \times (P \times \R)$ be a compact set and $W_s$ a family of proper exact Lagrangian submanifolds. Take any $p \in K$, $s\in[0,1]$, and sufficiently small $r>0$. There is a constant $D>0$ independent of these choices for which the following holds. Any non-constant $J_s$-holomorphic curve inside the ball $B_r(p)$ of radius $r$ which, moreover, satisfies the properties that
\begin{itemize}
\item it passes through the centre $p \in K$ of the ball, and
\item its boundary is located on $\partial B_r(p) \cup (B_r(p) \cap W_s)$,
\end{itemize}
has $\omega$-area bounded from below by $Dr^2$.

Using this monotonicity property, we will argue by contradiction, showing that if a strip as in (1), (2), and (3) leaves the set $U_0$, $U_{B+2}$, and $U_{[A-1,B]}$, respectively, then its $\omega$-area must exceed the corresponding upper bound found above.

To that end, there is a fixed metric on $\R \times (P \times \R)$ and a constant $C>0$ independent of $\epsilon$ and $\eta$ for which the following hold. The two components $(L \cup  L^{\eta,\epsilon}_\pm ) \cap \{ t = t_0 \}$, are at a distance at least $C\epsilon^2$ for $t_0 \in \{0,A-1,B\}$, while the two components $(L \cup  L^{\eta,\epsilon}_- ) \cap \{ t = B+2 \}$ are at a distance at least $C(\epsilon^{1/2}-\eta^2)$.

The monotonicity property for the $\omega$-area of $J_s$-holomorphic curves with boundary on $L \cup  L^{\eta,\epsilon}_\pm$ thus gives a constant $D>0$ independent of $\epsilon>\eta>0$ for which the following holds.
\begin{itemize}
\item If a $J_s$-holomorphic strip as above has a boundary-point passing through either of the sets
\begin{gather*}
V_0:=(L \cup  L^{\eta,\epsilon}_\pm) \cap \{ t=0 \},\\
V_{A-1}:=(L \cup  L^{\eta,\epsilon}_\pm) \cap \{ t =A-1 \},\\
V_B:=(L \cup  L^{\eta,\epsilon}_\pm) \cap \{ t =B \},
\end{gather*}
then it must have $\omega$-area bounded from below by $D\epsilon^4$.
\item If a $J_s$-holomorphic strip as above has a boundary-point passing through the set
\[V_{B+2}:=(L \cup  L^{\eta,\epsilon}_-) \cap \{ t =B+2 \},\]
then it must have $\omega$-area bounded from below by $D(\epsilon^{1/2}-\eta^2)^2$.
\end{itemize}
To that end, by examining the above bounds of distances, we may assume that the ball of radius $\epsilon^2/2$ centred at any point in $V_0 \cup V_{A-1} \cup V_B$ intersects $(L \cup  L^{\eta,\epsilon}_\pm)$ in a disc, and that the ball of radius $(\epsilon^{1/2}-\eta^2)/2$ centred at any point in $V_{B+2}$ intersects $L \cup  L^{\eta,\epsilon}_-$ in a disc.

After shrinking $\epsilon>\eta>0$, a comparison with the above upper bounds for the $\omega$-area of the strips of either type (1) and (3), we conclude that no such $J_s$-holomorphic strip has a boundary-point passing through either of the sets $V_0$, $V_{A-1}$, or $V_B$. Similarly, we may assume that no $J_s$-holomorphic strip of type (2) has a boundary-point passing through the set $V_{B+2}$.

(1): The above argument shows that the boundary of a $J_s$-holomorphic strip of this type is contained inside the set $\{ t \le 0 \}$. In particular, for each $\epsilon>\eta>0$ sufficiently small, its boundary may be assumed to be disjoint from some fixed $\delta$-neighbourhood of $\partial U_0$. Suppose that such a $J_s$-holomorphic strip is not contained entirely in $U_0$. Consequently, it has an interior point passing through $p \in \partial U_0$, where $p$ is of distance at least $\delta>0$ from the boundary.

(2): The above argument shows that the boundary of a $J_s$-holomorphic strip of this type is contained inside the set $\{ t \le B+2 \}$. In particular, for each $\epsilon>\eta>0$ sufficiently small, its boundary may be assumed to be disjoint from some fixed $\delta$-neighbourhood of $\partial U_{B+2}$. Suppose that such a $J_s$-holomorphic strip is not contained entirely in $U_{B+2}$. Consequently, it has an interior point passing through $p \in \partial U_{B+2}$, where $p$ is of distance at least $\delta>0$ from the boundary.

(3): The above argument shows that the boundary of a $J_s$-holomorphic strip of this type is contained inside the set $\{ A-1 \le t \le B \}$. In particular, for each $\epsilon>\eta>0$ sufficiently small, its boundary may be assumed to be disjoint from some fixed $\delta$-neighbourhood of $\partial U_{[A-1,B]}$. Suppose that such a $J_s$-holomorphic strip is not contained entirely in $U_{[A-1,B]}$. Consequently, it has an interior point passing through $p \in \partial U_{[A-1,B]}$, where $p$ is of distance at least $\delta>0$ from the boundary.

In all of the three cases above, the monotonicity property for the $\omega$-area of $J_s$-holomorphic curves without boundary shows that the strip must have $\omega$-area bounded from below by $D\delta^2>0$, where both $D>0$ and $\delta>0$ are independent of the choices of (sufficiently small) numbers $\epsilon>\eta>0$. When $\epsilon>\eta>0$ is small enough, this contradicts the above upper bounds of the $\omega$-area, which means that the $J_s$-holomorphic strip must be contained inside the required neighbourhood.
\end{proof}

\section{Proof of  Theorem \ref{thm:lift}}
\label{sec:proofs}

In the following, we fix a front-generic Legendrian submanifold $\Lambda \subset P \times \R$.  We begin with the below lemma, which is a generalisation of the analysis in \cite[Theorem 7.7]{InvLegCoh} to the case of an arbitrary contactisation. Here we make the identification $P \times \R \simeq P \times i\R \subset P \times \C$ and define
\[ L:=\{ (p,x+iy) \in P \times \C;\:\: (p,iy) \in \Lambda \} \simeq \Lambda \times \R.\]
\begin{lemma}
\label{lem:lift}
Given a $J_P$-holomorphic polygon $u \co (D^2,\partial D^2) \to (P,\Pi_{\mathrm{Lag}}(\Lambda))$ with punctures mapping to double-points, there is a $(J_P \oplus i)$-holomorphic lift
\[(u,a) \co (D^2,\partial D^2) \to (P\times \C,L).\]
Let $p \in \partial D^2$ be a puncture that is mapped by $u$ to a double-point corresponding to the Reeb chord $\{u(p)\} \times [A,A+\ell]$, and choose holomorphic coordinates identifying $D^2$ with $\{s+it; \: 0 \le t \le 1\} \subset \C$, in which $p$ corresponds to $s=+\infty$. It follows that there are $C \in \C$, $\lambda>0$, for which the bound
\[ \| a(s+it)-C \pm \ell(s+it) \| \le e^{-\lambda s}\]
holds for all $s > 0$ sufficiently large.
\end{lemma}
\begin{proof}
Away from the boundary punctures, there is a unique continuous lift
\[(u|_{\partial D^2},h) \co \partial D^2 \to \Lambda \subset P \times \R.\]
By abuse of notation, we let $h \co D^2 \to \R$ denote the harmonic extension to the interior. Observe that $h$ is bounded and $C^\infty$ away from the boundary punctures, where it has jump discontinuities.

Let $-g$ be the Harmonic conjugate of $h$, which is smooth in the interior of $D^2$. We define
\[a(s+it):=g(s,t)+ih(s,t),\]
which is holomorphic in the interior of $D^2$ and satisfies the correct boundary condition. It remains to show that $a$ has the above asymptotic behaviour at its boundary punctures.

In the following we will use $S_0 >0$ to denote some sufficiently large real number.

We fix a boundary puncture $p \in \partial D^2$ of $u$. Choose a conformal identification of the domain $D^2$ with the strip $\{ s+it; \: 0 \le t \le 1 \} \subset \C$ such that the boundary-puncture $p$ corresponds to $s=+\infty$. An application of \cite[Theorem B]{AsympStrips} yields that there are $C^\infty$ functions $\mathbf{v},\mathbf{w} \co \{ s+it; \:\: 0 \le t \le 1 \} \to \C^n$, and a number $\lambda>0$, such that
\[ u(s+it)= -\frac{1}{\lambda} e^{-\lambda s}\mathbf{v}(t)+\mathbf{w}(s,t), \:\: s\ge S_0,\]
holds in a Darboux chart centred at $u(p)$. Moreover, there is a $\delta>0$ such that
\[\|\mathbf{w}(s,t)\|_{C^k[s,+\infty)}= \mathbf{O}(e^{-(\lambda+\delta)s})\]
holds for $s \ge S_0$.

Let $h_+(s):=h(s,1)$ and $h_-(s)=h(s,0)$ denote $h$ above restricted to the respective boundary component of the strip. Recall that the contact form is given by $dz+\theta$, where $\theta$ is a one-form on $P$. Considering the expressions
\begin{gather*}
h_+(s_1)-h_+(s_0)=\int_{s_0}^{s_1} u|_{\{t=1\}}\,^*(-\theta) ds,\\
h_-(s_1)-h_-(s_0)=\int_{s_0}^{s_1} u|_{\{t=0\}}\,^*(-\theta) ds,
\end{gather*}
and using the above asymptotic expansion, it follows that
\begin{eqnarray}
h_\pm(s) & =&  C_\pm+\mathbf{O}(e^{-\lambda s}) \label{eq:h},\\
\left(\frac{d}{d s}\right)^{k}h_\pm(s) & = & \mathbf{O}(e^{-\lambda s}), \:\:k =1,2,3,\hdots \label{eq:hh}
\end{eqnarray}
holds for $s \ge S_0$.

The Poisson kernel for the strip $\{s+it; \:\: 0 \le t \le 1\}$ is given by
\[P(s,t)=\pi\frac{\sin{\pi t}}{\cosh{\pi s}-\cos{\pi t}},\]
as computed in \cite{HarmonicStrip}, and the harmonic extension $f(s,t)$ to the strip of a function given by $f_+(s)$ along $\{t =1\}$ and $f_-(s)$ along $\{t=0\}$ is given by the convolution
\begin{eqnarray}
\lefteqn{f(s,t)=} \label{eq:poisson} \\
& = & \frac{1}{2\pi} \int_{-\infty}^{+\infty} \left(P(\sigma-s,t)f_-(\sigma) +P(\sigma-s,1-t)f_+(\sigma) \right) d\sigma \nonumber \\
& = & \frac{1}{2\pi} \int_{-\infty}^{+\infty} \left(P(\sigma,t)f_-(\sigma+s) +P(\sigma,1-t)f_+(\sigma+s) \right) d\sigma. \nonumber
\end{eqnarray}

Observe that the harmonic extension $h(s,t)$ can be expressed by such a convolution and, by the regularity of $h_\pm$, it is $C^\infty$ up to the boundary for $s  \ge S_0$. Likewise, the bounded harmonic function
\[h_0(s,t):=C_-+(C_+-C_-)t\]
on the strip can expressed by a convergent convolution of the Poisson kernel with its restriction to the boundary. Since
\[h(s,t)-h_0(s,t)=\mathbf{O}(e^{-\lambda s}), \:\: t \in \{0,1\}, s \ge S_0, \]
holds along the boundary of the strip by (\ref{eq:h}), using formula (\ref{eq:poisson}) we get the estimate
\begin{equation}
h(s,t)-h_0(s,t)=\mathbf{O}(e^{-\lambda s}), \:\: s \ge S_0, \label{eq:hhh}
\end{equation}
in the interior of the strip as well.

Since convolving with the Poisson kernel commutes with the operation $\partial_s$, using (\ref{eq:hh}) it similarly follows that
\begin{equation}
\label{eq:as}
(\partial_s)^k h(s,t)=\mathbf{O}(e^{-\lambda s}),\:\: k > 0, s \ge S_0.
\end{equation}

In order to find a bound for $\partial_t h(s,t)$ we proceed as follows. Observe that
\begin{eqnarray*}
\lefteqn{\partial_t (h-h_0)(s,1/2) =}\\
& = & \frac{1}{2\pi} \int_{-\infty}^{+\infty} \partial_t P(\sigma,1/2)\left((h_-(\sigma+s)-C_-)-(h_+(\sigma+s)-C_+)\right) d\sigma,
\end{eqnarray*}
where
\[\partial_t P(s,1/2)=\pi^2\left.\frac{\cos{\pi t}\cosh{\pi s}-1}{(\cosh{\pi s}-\cos{\pi t})^2}\right|_{t=1/2}=\mathbf{O}(e^{-2|s|})\]
for $|s| \ge S_0$. Together with (\ref{eq:h}) we get the estimate
\begin{equation}
\label{eq:as2}
\partial_t (h-h_0)(s,1/2)=\mathbf{O}(e^{-\lambda s})
\end{equation}
for $s \ge S_0$. The harmonicity of $(h-h_0)(s,t)$ together with (\ref{eq:as}) implies that
\[\partial_t^2(h-h_0)(s,t)=-\partial_s^2(h-h_0)(s,t)=\mathbf{O}(e^{-\lambda s}),\]
which together with (\ref{eq:as2}) can be integrated to give
\[ \partial_t (h-h_0)(s,t)=\mathbf{O}(e^{-\lambda s}),\:\: s\ge S_0.\]
The latter identity is equivalent to
\begin{equation}
\partial_t h(s,t)=(C_+-C_-) + \mathbf{O}(e^{-\lambda s}), \label{eq:as3}
\end{equation}
for $s \ge S_0$. 

In particular, the Harmonic conjugate $-g(s,t)$ of $h(s,t)$ on the strip satisfies
\[g(s,t)=B+(C_+-C_-)s+\mathbf{O}(e^{-\lambda s}),\]
for some $B \in \R$, as follows by integrating the Cauchy-Riemann equations
\[\partial_s g =\partial_t h, \:\: \partial_t g=-\partial_s h,\]
together with the estimates (\ref{eq:as}) and (\ref{eq:as3}).

This estimate together with (\ref{eq:hhh}) finally gives the sought asymptotic bound of $a(s+it):=g(s,t)+ih(s,t)$, after possibly choosing a smaller $\lambda>0$.
\end{proof}

\begin{remark}
\label{rem:lift}
In the case when $\theta=-d^{J_P}\alpha =-d\alpha(J_P \cdot)$ for some smooth function $\alpha \co P \to \R$, the map
\begin{gather*}
(\R \times (P \times \R),\R \times \Lambda) \to (P \times \C,L) \\
(t,(p,z)) \mapsto (p,t-\alpha(p)+iz)
\end{gather*}
can be seen to pull back $J_P \oplus i$ to the cylindrical lift of $J_P$.
\end{remark}

In the following we fix a regular compatible almost complex structure $J_P$ on $P$ which moreover is integrable in some neighbourhood of the double-points of $\Pi_{\mathrm{Lag}}(\Lambda)$.

\begin{proof}[Proof of Theorem \ref{thm:lift}]
Since the projection $\pi_P$ is $(\widetilde{J}_P,J_P)$-holomorphic, and since $J_P$ is regular, Lemma \ref{lem:trans} below implies that $\widetilde{J}_P$ is regular as well and that the smooth map
\begin{gather*}
\mathcal{P} \co \mathcal{M}_{a;\mathbf{b}}(\R \times \Lambda;\widetilde{J}_P)/\R \to \mathcal{M}_{a;\mathbf{b}}(\Pi_{\mathrm{Lag}}(\Lambda);J_P), \\
\widetilde{u} \mapsto \pi_P \circ \widetilde{u}
\end{gather*}
is a local diffeomorphism. We will show that
\begin{enumerate}
\item $\mathcal{P}$ is injective, and
\item $\mathcal{P}$ takes values in every component of $\mathcal{M}_{a;\mathbf{b}}(\Pi_{\mathrm{Lag}}(\Lambda);J_P)$.
\end{enumerate}
Observe that, since $\mathcal{P}$ has a natural extension to a map between the respective compactifications of the above moduli spaces, it can be seen that the image of $\mathcal{P}$ is closed (we refer to the proof of part (2) below for more details about these compactifications). Together with the above properties, it follows that $\mathcal{P}$ is a diffeomorphism.

We begin by making the following observation. The integrability of $J_P$ in a neighbourhood of the double-points of $\Pi_{\mathrm{Lag}}(\Lambda)$, together with the compatibility of $J_P$, implies that $d\theta$ is a K\"{a}hler form there. It is well-known that $d\theta=-dd^{J_P}f$ holds in some possibly smaller neighbourhood, for some smooth real-valued function $f$. Locally we have $\theta=-d^{J_P}f+dg$ for some smooth real-valued function $g$. The function defined by
\[\alpha(\mathbf{z}):=f(\mathbf{z})+g(i \cdot \mathbf{z})\]
in some local holomorphic chart thus satisfies $-d^{J_P}\alpha=\theta$. This fact will be used in the proofs of both part (1) and (2).

(1): Since $\theta=-d^{J_P}\alpha$ holds in some neighbourhood $U \subset P$ of the double-points, the identification made in Remark \ref{rem:lift} can be applied to identify $(\R \times (U \times \R),\widetilde{J}_P)$ with $(U \times \C, J_P \oplus i)$.

Suppose that two $\widetilde{J}_P$-holomorphic discs $\widetilde{u}$ and $\widetilde{u}'$ satisfy
\[u:=\mathcal{P}(\widetilde{u})=\mathcal{P}(\widetilde{u}').\]
The boundary-condition implies that $\pi_\C \circ \widetilde{u}|_{u^{-1}(U)}$ and $\pi_\C \circ \widetilde{u}'|_{u^{-1}(U)}$ are holomorphic maps whose imaginary parts agree along the boundary, where we have used the above identification. By a standard result, these maps differ by a real constant on each component of $u^{-1}(U)$ that intersects $\partial D^2$.

Note that $u^{-1}(U) \cap \partial \dot{D}^2 \neq \emptyset$ by the behaviour of $u$ near the boundary-punctures. It thus follows that, after a translation of the $t$-coordinate, $\widetilde{u}$ and $\widetilde{u}'$ may be assumed to coincide on a non-empty open set. The union of critical points
\[P:=\Crit(\widetilde{u}) \cup \Crit(\widetilde{u}') \subset \dot{D}^2\]
is a discrete set by \cite[Theorem 3.5]{ExistenceInjective}. Now \cite[Lemma 4.2]{ExistenceInjective} can be used to show that the equality $\widetilde{u}=\widetilde{u}'$ holds on a both open and closed subset of $\mathrm{int}\dot{D^2} \setminus P$. In conclusion, $\widetilde{u}=\widetilde{u}'$ holds on $\dot{D}^2$, which implies the injectivity.

(2): In the case when the assumptions of Remark \ref{rem:lift} are satisfied, that is when $\theta=-d^{J_P}\alpha$, Lemma \ref{lem:lift} can be used to construct an explicit lift to a $J_P$-holomorphic disc, which immediately gives this property. In the general case we proceed as follows.

The above moduli spaces can be compactified by using Gromov-Hofer compactness. See \cite{ContHomP} for the case of the moduli space of $J_P$-holomorphic polygons in $P$ and \cite{CompSFT} for the case of the moduli space of $\widetilde{J}_P$-holomorphic discs in the symplectisation. The compactified moduli spaces are manifolds having boundary with corners.

Points in the boundary strata of these moduli spaces correspond to so called broken configurations. In this setting, a broken configuration can be identified with a connected directed tree with the extra data consisting of
\begin{itemize}
\item a pseudo-holomorphic disc in a moduli space as above assigned to each vertex,
\item a Reeb chord assigned to each directed edge,
\end{itemize}
and satisfying the following condition. Given an edge starting at the vertex $v_1$, ending at the vertex $v_2$, and corresponding to the Reeb chord $c$, we require the pseudo-holomorphic disc corresponding to $v_1$ (respectively $v_2$) to have a negative (respectively positive) boundary-puncture at $c$.

A broken configuration can be glued to give a pseudo-holomorphic disc. Recall that the sum of $d\theta$-energies of the discs in a broken configuration equals the $d\theta$-energies of the corresponding glued pseudo-holomorphic disc.

We let
\[\mathcal{P} \co \overline{\mathcal{M}}_{a;\mathbf{b}}(\R \times \Lambda;\widetilde{J}_P)/\R \to \overline{\mathcal{M}}_{a;\mathbf{b}}(\Pi_{\mathrm{Lag}}(\Lambda);J_P)\]
be the natural extension of $\mathcal{P}$ to a continuous map on the compactified moduli spaces, where the domain should be interpreted as the compactification of the quotient $\mathcal{M}_{a;\mathbf{b}}(\R \times \Lambda;\widetilde{J}_P)/\R$.

By the formula for the $d\theta$-energy in terms of the action of Reeb chords, together with the fact that there are only finitely many Reeb chords on $\Lambda$, it follows that there are finitely many non-empty moduli spaces as above in the case under consideration. Moreover, the above compactness theorems imply that each moduli space has finitely many components.

We prove the theorem by induction on the $d\theta$-energy. Assume that the statement has been shown for every moduli space consisting of discs having $d\theta$-energy strictly less than $A$. We must now show that $\mathcal{P}$ maps into each component of the moduli space $\mathcal{M}_{a;\mathbf{b}}(\Pi_{\mathrm{Lag}}(\Lambda);J_P)$ consisting of $J_P$-holomorphic polygons having $d\theta$-energy equal to $A$.

We begin with the case of a non-compact component of the moduli space $\mathcal{M}_{a;\mathbf{b}}(\Pi_{\mathrm{Lag}}(\Lambda);J_P)$, i.e.~a component whose compactification has non-empty boundary. Since each pseudo-holomorphic disc in a broken configuration corresponding to a boundary-point has $d\theta$-area strictly less than $A$, the induction hypothesis implies that the above map $\mathcal{P}$ maps into its boundary. Gluing a broken configuration in $\overline{\mathcal{M}}_{a;\mathbf{b}}(\R \times \Lambda;\widetilde{J}_P)/\R$ being a point in the preimage of such a boundary-point produces a solution contained in the interior of $\mathcal{M}_{a;\mathbf{b}}(\R \times \Lambda;\widetilde{J}_P)/\R$. Furthermore, this solution projects to the corresponding component of $\mathcal{M}_{a;\mathbf{b}}(\Pi_{\mathrm{Lag}}(\Lambda);J_P)$, which finishes the claim in this case.

The remainder of the proof concerns the case of a compact component of $\mathcal{M}_{a;\mathbf{b}}(\Pi_{\mathrm{Lag}}(\Lambda);J_P)$, i.e.~a component that is a closed manifold. In this case we proceed as follows to show the existence of a lifted solution.

Recall the existence of the smooth function $\alpha \co U \to \R$ satisfying $\theta=-d^{J_P}\alpha$, where $U \subset P$ is a neighbourhood of the double-points. We choose a cut-off function $\rho \co P \to \R_{\ge 0}$ having support inside $U$ and which has the property that $\rho=1$ in a neighbourhood $V \subset U$ of the double-points. We construct the smooth function
\begin{gather*}
\beta \co P \times \R \to \R,\\
\beta(p):=\begin{cases} \rho(p)\cdot \alpha(p), & p \in U \\
0, & p \notin U.
\end{cases}
\end{gather*}

Consider the one-parameter family of smooth one-forms $\lambda_s$ on $P \times \R$ defined by
\[ \lambda_s(p,z) := \begin{cases} dz+s\theta, & p \in P \setminus U,\\
dz-d^{J_P}((1-s)\beta+s\alpha), & p \in U.
\end{cases}\]
Observe that $\lambda_1=dz+\theta$ coincides with the original contact form on $P \times \R$ and that each $\lambda_s$ coincides with the original contact form in the neighbourhood $V \times \R$ containing the Reeb chords.

The induced one-parameter family $\xi_s:=\ker \lambda_s \subset T(P \times \R)$ of tangent hyper-plane fields satisfies the following properties for each $s \in [0,1]$.
\begin{itemize}
\item $\xi_s$ is transverse to $\partial_z$.
\item $\xi_s$ is invariant with respect to translations of the $z$-coordinate.
\item $\xi_s$ coincides with the contact-distribution $\ker (dz+\theta)$ in $V \times \R$.
\item $d\theta$ restricts to a symplectic form on $\xi_s$.
\end{itemize}

For every $s \in [0,1]$ we lift $J_P$ to a an almost complex structure $J_s$ on the symplectisation $\R \times (P \times \R)$ which is uniquely determined by the requirements that
\begin{itemize}
\item $J_s$ is invariant with respect to translations of the $t$ and $z$-coordinates,
\item $J_s \partial_t = \partial_z$,
\item $J_s \xi_s = \xi_s$, and
\item $\pi_P$ is $(J_s,J_P)$-holomorphic.
\end{itemize}

First observe that, since $\xi_1=\ker(dz+\theta)$ is the contact distribution, by construction we have $J_1=\widetilde{J}_P$. For the same reasons, since $\xi_s$ agrees with the contact distribution in $V \times \R$, it follows that $J_s=\widetilde{J}_P$ in $\R \times (V \times \R)$ for each $s \in [0,1]$.

Second, it can be checked that $J_0$ coincides with the pull-back of the almost complex structure $J_P \oplus i$ on $P \times \C$ under the diffeomorphism
\begin{gather*}
\iota\co\R \times (P \times \R) \to P \times \C, \\
(t,(p,z)) \mapsto (p,t-\beta(p)+iz).
\end{gather*}
To that end, observe that the hypersurfaces $\{t \} \times P \times \R$ are mapped by $\iota$ to level-sets of the smooth function 
\begin{gather*}
\tau \co P \times \C \to \R,\\
(p,x+iy) \mapsto x+\beta(p),
\end{gather*}
and that, since
\[ \iota^*( -d^{J_P \oplus i}\tau )=\iota^*( dy-d^{J_P}\beta )=dz-d^{J_P}\beta=\lambda_0,\]
the $(J_P \oplus i)$-complex tangencies to the latter hypersurfaces correspond to $J_0$-complex tangencies to the former.

Recall that we want to lift a $J_P$-holomorphic disc lying in a transversely cut out closed component $\mathcal{M} \subset \mathcal{M}_{a;\mathbf{b}}(\Pi_{\mathrm{Lag}}(\Lambda);J_P)$. By using Lemma \ref{lem:lift} together with the above $(J_0,J_P \oplus i)$-holomorphic diffeomorphism $\iota$ we can lift each $u \in \mathcal{M}$ to a finite-energy $J_0$-holomorphic disc in the symplectisation. We use $\widetilde{\mathcal{M}} \subset \mathcal{M}_{a;\mathbf{b}}( \R \times \Lambda;J_0)$ to denote the component of the moduli space containing these lifts, which is transversely cut out by Lemma \ref{lem:trans}. Observe that $\widetilde{\mathcal{M}}$ is a closed manifold and that the above projection induces a diffeomorphism $\widetilde{\mathcal{M}}  \to \mathcal{M}$.

Consider the connected component 
\[\mathcal{W} \subset \bigcup_{s \in [0,1]} \mathcal{M}_{a;\mathbf{b}}( \R \times \Lambda;J_s)\]
containing $\widetilde{\mathcal{M}}$. Again, the above projection induces a map $\widetilde{\mathcal{P}} \co \mathcal{W}  \to \mathcal{M}$,
which thus is of degree one when restricted to $\mathcal{W} \cap \{s=0 \}=\widetilde{\mathcal{M}}$.

Lemma \ref{lem:trans} implies that $\mathcal{W}$ is transversely cut out, and is hence a smooth finite-dimensional manifold. The goal is to show that $\mathcal{W}$ is compact. In this case, since $\mathcal{W}$ is a compact cobordism from the closed manifold $\mathcal{W} \cap \{s=0\}=\widetilde{\mathcal{M}}$ to the closed manifold $\mathcal{W} \cap \{s=1\}$, it follows that $\widetilde{\mathcal{P}}$ restricted to $\mathcal{W} \cap \{s=1\}$ has degree one as well. In particular, $\mathcal{W} \cap \{s=1\}$ is non-empty, and there is a $\widetilde{J}_P$-holomorphic lift of any solution inside $\mathcal{M}$.

To show the compactness of $\mathcal{W}$ we argue as follows. First, recall that $J_s=\widetilde{J}_P$ in a neighbourhood of the form $\R \times (V \times \R)$ where $V \times \R$ contains all the Reeb chords on $\Lambda$. In particular, outside of some compact subset of the domain, two solutions in $\mathcal{W}$ are solutions to the same boundary-value problem.

Second, since we have
\[\R \partial_z = \ker d\theta \subset \ker d\lambda_s,\:\: \iota_{\partial_z}\lambda_s=1,\]
the pair $(d\theta,\lambda_s)$ is a so-called stable Hamiltonian structure on $P \times \R$. Moreover, $J_s$ is compatible with this stable Hamiltonian structure in the sense described in \cite[Section 2.2]{CompSFT}. In conclusion, it follows that the Gromov-Hofer compactness applies to the moduli-space $\mathcal{W}$ of $J_s$-holomorphic curves.

Furthermore, the computation
\[ -d^{J_s}(t-(1-s)\beta)=\lambda_s+(1-s)d^{J_P}\beta=dz+s\theta\]
shows that the function $t-(1-s)\beta$ on the symplectisation is weakly $J_s$-convex. This implies that the maximum-principle holds for this function pulled back to any $J_s$-holomorphic curve, which prevents certain bad breakings.

The above shows that $\mathcal{W}$ may be compactified as usual. However, no sequence of curves in $\mathcal{W}$ can converge to a broken solution, since such a solution necessarily would project under $\widetilde{\mathcal{P}}$ to a broken solution in the boundary of $\overline{\mathcal{M}} \subset \mathcal{M}_{a;\mathbf{b}}(\Pi_{\mathrm{Lag}}(\Lambda);J_P)$, which is empty by assumption. It thus follows that $\mathcal{W}$ was compact to begin with.
\end{proof}

\section{Transversality for lifted discs}
Let $\pi_P$ be $(J,J_P)$-holomorphic. In this section we show that the transversality of certain moduli spaces of $J$-holomorphic discs $\widetilde{u}$ in the symplectisation of $P \times \R$ follows from the transversality of the corresponding moduli space of $J_P$-holomorphic discs
\[\mathcal{P}(\widetilde{u}):=\pi_P\circ \widetilde{u}. \]

It should be mentioned that we do not specify which functional-analytic set-up we are using. On the other hand, the proofs use arguments applied to solutions of a linearised problem which are smooth, as follows by elliptic regularity. Hence, these arguments can be applied independently of the functional-analytic set-up.

\subsection{A sketch of the construction of the linearised operator}
Here we sketch the construction of the linearisations of the non-linear Cauchy-Riemann operators under consideration. 

The moduli spaces of the above $J_P$-holomorphic discs are constructed as the vanishing-locus of the section
\[u \mapsto \overline{\partial}_{J_P}u=du+J_P du \circ i\]
of a bundle $\mathcal{E}^{0,1} \to \mathcal{B}$. Here $\mathcal{B}$ denotes a Banach manifold consisting of maps $u \co (D^2,\partial D^2) \to (P,\Pi_{\mathrm{Lag}}(\Lambda))$ in an appropriate Sobolev class. We refer to \cite{ContHomP} for more details.

By the linearisation of $\overline{\partial}_{J_P}$ at a $J_P$-holomorphic disc $u \in \mathcal{B}$ we denote the linear operator
\[ D_u := \pi_\mathcal{F} \circ T\overline{\partial}_{J_P} \co  T_u\mathcal{B} \to \Omega^{0,1}(u^*TP), \]
where $\pi_\mathcal{F}$ is the projection to the fibre over $u$ of the bundle $\mathcal{E}^{0,1}$.

By the choice of a metric, one can make the identification 
\[T_u\mathcal{B} \simeq \Gamma\left(u^*TP,u|_{\partial D^2}^*T(\Pi_{\mathrm{Lag}}(\Lambda))\right),\]
 where the right-hand side denotes the space of sections in $u^*TP$ which take value in the sub-bundle $u|_{\partial D^2}^*T(\Pi_{\mathrm{Lag}}(\Lambda))$ along the boundary. Observe that we have to consider the completion of the space of smooth sections with respect to an appropriate Sobolev norm, but that we have suppressed this from the notation.

In the the case of the symplectisation $\R \times (P \times \R)$ there are similar constructions, where $\overline{\partial}_J$ becomes a section in the bundle $\widetilde{\mathcal{E}}^{0,1} \to \widetilde{\mathcal{B}}$ over the Banach manifold $\widetilde{\mathcal{B}}$ consisting of maps $\widetilde{u} \co (\dot{D}^2,\partial \dot{D}^2) \to (\R \times (P \times \R),\R \times \Lambda )$ in an appropriate Sobolev class (with positive weights). We refer to \cite{StripsII} and \cite{FredholmTheory} for examples of this, although the analytical set-up there differs from that in \cite{ContHomP}.

In any case, when $\mathcal{B}$ and $\widetilde{\mathcal{B}}$ consists of maps in suitable Sobolev spaces, the linearisation $D_u$ at a solution $u$ of the respective non-linear Cauchy-Riemann operator is elliptic. We use
\[\indx D_u = \dim\ker D_u- \dim\coker D_u\]
to denote its Fredholm index. In the case when the linearisation $D_u$ is surjective, that is $\coker D_u=0$, it follows that a neighbourhood of $u$ in its moduli space is transversely cut out. In particular, this neighbourhood is a smooth finite-dimensional manifold.

\subsection{The transversality results}

Let $J$ be a (not necessarily cylindrical) almost complex structure on $\R \times (P \times \R)$ that is invariant under translation in the $t$ and $z$-coordinates, satisfies $J \partial_t = \partial_z$, and which has the property that
\[ \pi_P \co \R \times (P \times \R) \to P\]
is $(J,J_P)$-holomorphic. Furthermore, we assume that $J$ preserves the contact-planes above neighbourhoods of the double-points of $\Pi_{\mathrm{Lag}}(\Lambda) \subset P$. 

We fix the metric
\begin{gather*}
\widetilde{g}:=dt \otimes dt+(\pi_P)^*(g)+(dz+\theta)\otimes (dz+\theta), \\
g(v_1,v_2) := d\theta(v_1,J_Pv_2), \:\: v_i \in TP,
\end{gather*}
on $\R \times (P \times \R)$, where $g$ is the metric on $P$ induced by the symplectic form and choice of compatible almost complex structure. Observe that under the splitting
\[T_{(t_0,(p_0,z_0))}(\R \times (P \times \R))=\R \partial_t \oplus \ker (dz+\theta) \oplus \R \partial_z \]
the exponential map induced by $\widetilde{g}$ takes the form
\begin{gather*}
\widetilde{\exp}_{(t_0,(p_0,z_0))}(t,v,0)=(t_0+t,(\exp_{p_0}(v),z_0+F(p_0,v))), \\
\widetilde{\exp}_{(t_0,(p_0,z_0))}(t,0,z)=(t_0+t,p_0,z_0+z),
\end{gather*}
where $\exp$ denotes the exponential map on $P$ induced by the metric $g$, and where $F \co P \times \ker (dz+\theta) \to \R$ is a smooth function.

We use $\widetilde{\exp}$ to make the identification
\[T_{\widetilde{u}}\widetilde{\mathcal{B}} \simeq \Gamma\left(\widetilde{u}^*T(\R \times (P \times \R)),\widetilde{u}|_{\partial D^2}^*(T(\R \times \Lambda))\right)\]
 and $\exp$ to make the identification
\[T_u\mathcal{B} \simeq \Gamma\left( u^*TP,u|_{\partial D^2}^*(T\Pi_{\mathrm{Lag}}(\Lambda))\right).\]
Writing $u:=\pi_P \circ \widetilde{u}$, the tangent-map
\[T \pi_P \co \Gamma(\widetilde{u}^*T(\R \times (P \times \R))) \to \Gamma(u^*TP)\]
has a canonical extension to
\[T \pi_P \co \Omega^{0,1}(\widetilde{u}^*T(\R \times (P \times \R))) \to \Omega^{0,1}(u^*TP).\]
Using the properties
\begin{gather*}
\pi_P \circ \widetilde{\exp}(0,v,0)=\exp(v),\\
\pi_P \circ \widetilde{\exp}(t,0,z)=\exp(0),\\
T\pi_P (\overline{\partial}_J \widetilde{u})=T\pi_P(d\widetilde{u}+Jd\widetilde{u}\circ i)=\overline{\partial}_{J_P}(\pi_P \circ \widetilde{u})
\end{gather*}
it readily follows that
\begin{lemma}
\label{lem:tp}
For the above choices, and $u=\pi_P \circ \widetilde{u}$, we have
\[T \pi_P \circ D_{\widetilde{u}} = D_u \circ T \pi_P.\]
In particular, the restriction of $T\pi_P$ induces a linear map
\[ T\mathcal{P} \co \ker D_{\widetilde{u}} \to \ker D_u.\]
\end{lemma}

\begin{lemma}
\label{lem:trans}
Let $J$ be as above, and consider a finite-energy $J$-holomorphic disc
\[ \widetilde{u} \co(\dot{D}^2,\partial \dot{D}^2) \to (\R \times (P \times \R),\R \times \Lambda)\]
with punctures asymptotic to Reeb chords, and write $u:=\pi_P\circ \widetilde{u}$. It follows that
\[\indx {D_{\widetilde{u}}}=\indx{D_u}+1.\]
If $D_u$ is surjective it follows that $D_{\widetilde{u}}$ is surjective as well and, moreover, that the tangent-map
\[ T \mathcal{P} \co T_{\widetilde{u}}( \mathcal{M}_{a;\mathbf{b}}(\R \times \Lambda;J)/\R) \to T_u(\mathcal{M}_{a;\mathbf{b}}(\Pi_{\mathrm{Lag}}(\Lambda);J_P))\]
induced by the projection $\mathcal{P}$ is an isomorphism.
\end{lemma}
\begin{proof}
The statement concerning the indices follow from their expressions in terms of the gradings of the Reeb chords (see Section \ref{sec:dim}).

We now investigate $\ker T\mathcal{P}=\ker D_{\widetilde{u}} \cap \ker T \pi_P$ where
\[\ker T \pi_P \subset \Gamma\left(\widetilde{u}^*T(\R \times (P \times \R)),\widetilde{u}|_{\partial D^2}^*(T(\R \times \Lambda))\right)\]
consists of sections in the one-dimensional trivial complex sub-bundle whose fibre is given by $\R \langle \partial_t,\partial_z \rangle$. 

Using the fact that $J$ is invariant under translation of the $t$ and $z$-coordinates, it follows that $D_{\widetilde{u}}|_{\ker T \pi_P} $ can be identified with the standard Cauchy-Riemann operator on this trivial complex bundle. Hence, $\ker T\mathcal{P} \subset \ker T \pi_P$ can be seen to consist of holomorphic sections.

Using the identification of $x+iy \in \C$ with $x\partial_t+y\partial_z$, such a section is moreover identified with a holomorphic function $h \co \dot{D}^2 \to \C$. By elliptic regularity for the solutions in $\ker D_{\widetilde{u}}$, we may suppose that $h$ extends to a continuous function on the closed unit disc.

Since the linearised boundary condition moreover implies that $h$ is real (i.e.~takes values in $\R\partial_t$) when restricted to the boundary, it now follows that $h$ is constant and real. In other words, 
\[\ker T\mathcal{P} = \R \partial_t\]
is one-dimensional.

In conclusion, we have shown that
\[ 1+\dim \ker D_u \ge \dim \ker D_{\widetilde{u}}.\]
In the case when $\coker D_u=0$, we have
\[ \indx D_{\widetilde{u}} = 1+\indx D_u = 1+\dim \ker D_u,\]
which together with the above inequality yields $\indx D_{\widetilde{u}} \ge \dim \ker D_{\widetilde{u}}$. Finally, this implies that $\indx D_{\widetilde{u}} = \dim \ker D_{\widetilde{u}}$, and thus that $\coker D_{\widetilde{u}}=0$. From this it also follows that $T\mathcal{P}$ is an isomorphism.
\end{proof}

Let $V$ be an exact Lagrangian cobordism arising as $\phi^1_H(\R \times \Lambda')$, where the Hamiltonian flow
\[ \phi^s_H(t,(p,z))=(t,(p,z-s\sigma(t)))\]
is as described in Section \ref{sec:Reebflow}. In other words, $\sigma(t),\sigma(t)' \ge 0$, $\sigma(t)$ has support contained in $[N,+\infty)$, and $\sigma'(t)$ has compact support. We also assume that $(\R \times \Lambda) \cup V$ only has transverse double-points and that its Legendrian ends are chord-generic. For the cylindrical lift $\widetilde{J}_P$ to $\R \times (P \times \R)$ of a compatible almost complex structure $J_P$ on $P$, observe that the projection again induces a map
\begin{gather*}
\mathcal{P} \co \mathcal{M}_{f;\mathbf{a},e,\mathbf{b}}( (\R \times \Lambda) \cup V;\widetilde{J}_P) \to \mathcal{M}_{f;\mathbf{a},e,\mathbf{b}}( \Pi_{\mathrm{Lag}}(\Lambda \cup \Lambda');J_P), \\
\widetilde{u} \mapsto \pi_P \circ \widetilde{u},
\end{gather*}
of moduli spaces, where $f$ and $e$ denote either Reeb-chords or double-points.

\begin{lemma}
\label{lem:trans2}
Let $\widetilde{J}_P$ be the cylindrical lift of $J_P$ and consider $\widetilde{u} \in \mathcal{M}_{f;\mathbf{a},e,\mathbf{b}}( (\R \times \Lambda) \cup V;\widetilde{J}_P)$ and $u:=\pi_P\circ\widetilde{u}$ as above. It follows that
\[\indx {D_{\widetilde{u}}}=\indx{D_u}\]
in the case when $e$ is a double-point and
\[\indx {D_{\widetilde{u}}}=\indx{D_u}+1\]
in the case when $e$ is a Reeb chord. Moreover, the following holds:
\begin{enumerate}
\item Suppose that $J_P$ is integrable in a neighbourhood of the double-points of $\Pi_{\mathrm{Lag}}(P)$. If $\pi_P \circ \widetilde{u}$ is constant, i.e.~$\widetilde{u}$ is a strip contained a plane of the form $\R \times (\{q\} \times \R)$, it follows that $\widetilde{u}$ is transversely cut out.
\item If the negative puncture $e$ is mapped to a double-point of $(\R \times \Lambda) \cup V$ and if $J_P$ is regular, then $\widetilde{u}$ is transversely cut out as well.
\end{enumerate}
\end{lemma}

\begin{proof}
The statement concerning the indices again follow from their expressions in terms of the gradings of the Reeb chords and the double-points (see Sections \ref{sec:dim} and \ref{sec:transfer}).

Observe that Lemma \ref{lem:tp} applies in this setting as well, that is, the projection $T\pi_P$ induces a map
\[ T\mathcal{P} \co \ker D_{\widetilde{u}} \to \ker D_u.\]

(1): We assume that the domain of $\widetilde{u}$ is $D^2\setminus \{p_0,p_1\}$ and that its image is contained in the plane $\R \times (\{q\} \times \R)$ living over a double-point $q \in \Pi_{\mathrm{Lag}}(\Lambda \cup \Lambda') \subset P$. We use $L_\pm \subset T_qP$ denote the tangent-planes of the two branches of $\Pi_{\mathrm{Lag}}(\Lambda \cup \Lambda')$ at $q$. Furthermore, we make a conformal identification of the domain $D^2\setminus \{p_0,p_1\}$ with $\{s+it; \:\: t \in[0,1]\} \subset \C$.

We consider an element $\widetilde{\zeta}(s+it) \in \ker D_{\widetilde{u}}$ and let $\zeta(s+it):=T\mathcal{P}(\widetilde{\zeta})$. We begin by showing that $\zeta(s+it)=0$.

By the integrability assumption, there is a holomorphic chart near $q \in (P,J_P)$ which moreover induces an identification $T_qP \simeq \C^n$. Without loss of generality we may assume that the subspace $L_-$ is identified with the real-part $\R^n \subset \C^n$ while $L_+$ is the real span of $\{e^{i \theta_j}\mathbf{e}_j\}_{j=1,\hdots,n}$, where $\mathbf{e}_i$ denotes the standard basis of $\C^n$, and $0 < \theta_i < \pi$.

Because of the integrability of $J_P$ near $q$, the linearisation $D_u$ at the constant solution $u \co D^2 \setminus \{p_0,p_1\} \to \{q\} \subset P$ is the standard Cauchy-Riemann operator on sections $\zeta \co D^2 \setminus \{p_0,p_1\} \to \C^n \simeq T_qP$. In the above choice of holomorphic coordinates, it follows that the $j$:th component $\zeta_j(s+it)$ solves the boundary-value problem
\[\label{eq:lin}\begin{cases}
\partial_s \zeta_j(s+it) +i\partial_t\zeta_j(s+it)=0, \\
\zeta_j(s+i0) \in \R, \\
\zeta_j(s+i1) \in \R e^{i \theta_j}.
\end{cases}\]
In particular $\zeta_j(s+it)$ is holomorphic. Again, elliptic regularity implies that $\zeta$ has a continuous extension to the boundary of the unit disc. Since $\theta_j \neq 0$, the open mapping theorem shows that the holomorphic function $\zeta_j$ satisfying the above boundary conditions must vanish identically.

In conclusion, an element $\widetilde{\zeta} \in D_{\widetilde{u}}$ satisfies
\[T\mathcal{P}(\widetilde{\zeta})=0\]
and must hence be tangent to $\widetilde{u}$. In other words, it corresponds to an infinitesimal conformal reparametrisation of the domain, from which the claim follows.

(2): As in the proof of Lemma \ref{lem:trans} we will investigate $\ker T\mathcal{P}=\ker D_{\widetilde{u}} \cap \ker T \pi_P$. Again, it follows that $\ker T \pi_P \subset \Gamma (\widetilde{u}^*T(\R \times (P \times \R)))$ consists of the sections of the trivial complex bundle whose fibre $x\partial_t+y\partial_z$ can be identified with $x+iy \in \C$. Recall that, since $\widetilde{J}_P$ is invariant under translation of the $t$ and $z$-coordinates, elements of $\ker T\mathcal{P}$ can be identified with holomorphic functions $h \co \dot{D^2} \to \C$ that, by elliptic regularity, have a continuous extension to the boundary of the closed unit disc.

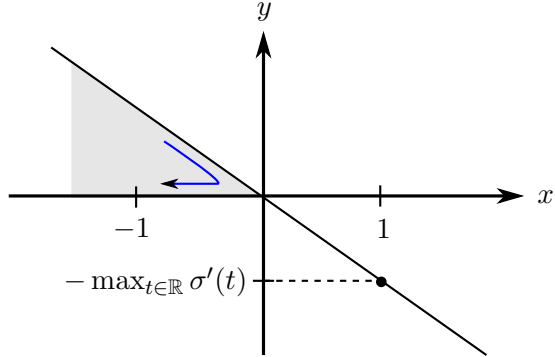
\begin{figure}[htp]
\centering
\labellist
\pinlabel $y$ at 105 130
\pinlabel $x$ at 211 60
\pinlabel $-1$ at 55 48
\pinlabel $1$ at 150 48
\pinlabel $-\max_{t\in\R}\sigma'(t)$ at 64 28
\endlabellist
\includegraphics{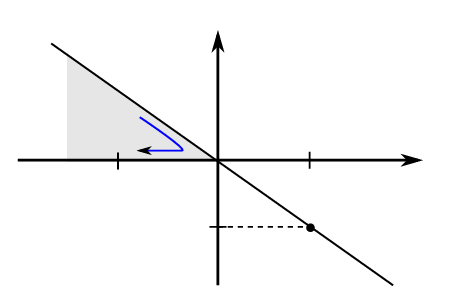}
\caption{The shaded region is the cone $(-1)^k\R_{\ge 0}(1-i\sigma'(\R))$ for $k=1$. The arrow schematically depicts the behaviour of $h$ along the (oriented) boundary near the negative puncture $p$ when $k=1$.}
\label{fig:corner}
\end{figure}

The linearised boundary condition implies that
\[ \begin{cases}
h(s+i0) \in \R, \\
h(s+i1) \in \R(1-i\sigma'(\R)),
\end{cases}\]
where we have used an appropriate identification of $D^2 \setminus \{ p_0,p_1\}$ with $\{s+it; \:\:t \in[0,1]\} \subset \C$. By using the open-mapping theorem it can be seen that $h$ actually maps into the cone $(-1)^k\R_{\ge 0}(1-i\sigma'(\R))$ for some $k \in \{0,1\}$ (recall that $\sigma'(\R) \subset \R_{\ge 0}$).

The linearised boundary condition has a discontinuity at the negative boundary-puncture $p \in \partial D^2$ mapping to the double-point $e$. It follows that $h$ must vanish there. Suppose that the $t$-coordinate takes the value $t_0$ at the double-point $e$ and observe that $\sigma'(t_0)>0$ by construction. Examining the linearised boundary condition near the puncture $p$ one sees that, for each $\sigma'(t_0)>\delta>0$, there is some $\epsilon>0$ for which
\[h(e^{i\theta}p) \in \begin{cases}
(-1)^k \R_{\ge 0}(1-i(\sigma'(t_0)+(-\delta,\delta))), & 0 \ge \theta>-\epsilon,\\
(-1)^k \R_{\ge 0}, & \epsilon> \theta \ge 0,
\end{cases}\]
where we have holomorphically identified the domain of $h$ with the unit disc. See Figure \ref{fig:corner} for a schematic picture.

Unless $h$ vanishes identically, a Schwarz-reflection together with the open mapping theorem shows that the image of $h$ must intersect the open quadrant
\[\{ x+iy;  \:\:x,y \in (-1)^k\R_{> 0}\}.\]
Together with $\sigma'(\R) \subset \R_{\ge 0}$, this contradicts the fact that the image of $h$ is contained inside of the cone $(-1)^k\R_{\ge 0}(1-i\sigma'(\R))$.

In conclusion, we have shown that $T\mathcal{P}$ is injective and hence that
\[ \indx D_{\widetilde{u}}=\indx D_u=\dim \ker D_u \ge \dim  \ker D_{\widetilde{u}}.\]
This implies that $\dim \coker D_{\widetilde{u}}=0$, i.e.~that $\widetilde{u}$ is transversely cut out.
\end{proof}

\begin{remark} The integrability condition in part (1) of the above lemma can be dropped. Using a separation of variables, the corresponding linearised equation can be solved in the general case as well. By analysing such a solution it can be shown that its projection $\zeta(s+it)$ either vanishes identically or blows up near a puncture. In the above proof, the open mapping theorem is used in order to deduce this.
\end{remark}

\appendix
\section{}

The following standard lemma is used in the invariance proof.
\begin{lemma}
\label{lem:hamiltonian}
Let $V$ be an exact Lagrangian cobordism from $\Lambda_-$ to $\Lambda$ in the symplectisation $(\R \times Y,d(e^t\lambda))$. Let $V'=\phi^1_{H_s}(V)$ be an exact Lagrangian cobordism from $\Lambda_-$ to $\Lambda'$ which is Hamiltonian isotopic to $V$ (we allow $\Lambda_- = \emptyset$).

If $\phi^s_{H_s}$ is supported in $[-N,+\infty) \times K$ for some compact set $K \subset Y$ and $N>0$, then there exists an exact Lagrangian cobordism $W$ from $\Lambda$ to $\Lambda'$ satisfying that
\begin{itemize}
\item $V \odot W$ is isotopic to $V'$ by a compactly supported Hamiltonian isotopy, and
\item $W$ is Hamiltonian isotopic to $\R \times \Lambda$ by a Hamiltonian isotopy as above.
\end{itemize}
\end{lemma}
\begin{remark} In particular, this lemma can be applied to the Hamiltonian isotopy of $W$ that it produces.
\end{remark}
\begin{proof}
Take a sufficiently large number $A > 0$ such that $V \cap \{ t \ge A \}$ and $V'  \cap \{ t \ge A \}$ both are cylindrical. Consider a smooth cut-off function $\rho \co \R \to \R$
 satisfying $0 \leq \rho(t) \leq 1$, $\rho|_{(-\infty,A]}=1$, and whose support is contained inside $(-\infty,B]$ for some $B>A$.

We define
\[W_s := ((-\infty,A] \times \Lambda) \cup \phi^s_{(1-\rho)H}(V \cap \{ t \ge A \}),\]
and observe that $W_0=\R \times \Lambda$, while $W:=W_1$ is an exact Lagrangian cobordism from $\Lambda$ to $\Lambda'$.

To infer that $W_1$ is cylindrical outside of a compact set we have used the assumption that $H_s$ has support inside $\R \times K$ where $K$ is compact, which implies that $\phi^s_H([-k,k] \times K) \subset \R \times K$, $k=1,2,3,\hdots$ is an exhaustion by compact sets.

We can make the identification
\[V \odot W = \phi^1_{(1-\rho)H}(V),\]
where the right-hand side obviously is exact Lagrangian isotopic to $V'$ via the compactly supported isotopy
\[V_\lambda:=\phi^1_{(1-\lambda\rho)H}(V), \:\: V_0=V', \: V_1=V \odot W \]
parametrised by $\lambda \in [0,1]$.

Finally, $W$ is Hamiltonian isotopic to $\R \times \Lambda$ via a Hamiltonian isotopy having support in $[A,+\infty) \times K$ by construction.
\end{proof}

For the following lemma, we assume that $g$ is a metric on $\Lambda$ for which $(f,g)$ constitutes a Morse-Smale pair. We let $(C^\bullet_{\mathrm{Morse}}(f),d_f)$ denote the induced Morse co-complex. Recall that $L \subset \R \times (P \times \R)$ is a manifold with a cylindrical end
\[L \cap \{ t \ge -2 \} = [-2,+\infty) \times \Lambda \subset \R \times (P \times \R),\]
and that $F_\pm$ are Morse-functions on $L$ as defined in Section \ref{sec:morse}.

\begin{lemma}
\label{lem:morse}
Consider a metric on $L$ which coincides with $dt \otimes dt + g$ on $\{ t \ge 0\}$. After a perturbation of this metric and of the function $F_-$, we may assume that the induced Morse co-complex satisfies
\[ (C^\bullet_{\mathrm{Morse}}(F_-),d_{F_-})=(C^\bullet_{\mathrm{Morse}}(F_+) \oplus C^{\bullet-1}_{\mathrm{Morse}}(f),d_{F_-}),\]
where
\[d_{F_-}= \begin{pmatrix} d_{F_+} & 0 \\
\Gamma & d_f
\end{pmatrix}.\]
In other words, $(C^\bullet_{\mathrm{Morse}}(F_-),d_{F_-})=\Cone(\Gamma)$ is the mapping-cone of a chain-map
\[ \Gamma \co (C^\bullet_{\mathrm{Morse}}(F_-),d_{F_-}) \to (C^\bullet_{\mathrm{Morse}}(f),d_f).\]
\end{lemma}
\begin{proof}
By considering the action-filtration of this Morse co-complex it follows that the sub-space
\[C^{\bullet-1}_{\mathrm{Morse}}(f) \subset C^\bullet_{\mathrm{Morse}}(F_-)\]
is a sub-complex. The only non-trivial part of the above statement is thus that the differential $d_{F_-}$ restricted to this sub-complex
may be assumed to coincide with $d_f$.

We restrict our attention to the domain
\[U:=[0,B] \times \Lambda \subset L \subset \R \times (P \times \R),\]
on which we let $s:=e^t|_U$ be the restriction. We construct the one-parameter family 
\[F_\lambda(s,p):=\epsilon^2 \alpha(s)+((1-\lambda)s+\lambda\alpha(s)) \eta^2 f(p), \:\:\lambda \in [0,1], \]
of smooth functions on $U$, where $\alpha(s)$ and $f$ are as constructed in Section \ref{sec:morse}. One can check that $F_\lambda$ is a one-parameter family of functions all whose negative gradients point outwards of $U$, and such that $F_0=F_-|_U$.

Recall that $\alpha(s)=s$ for $s \le e^{A-1}$, $\alpha(e^{B+1/2})=0$, and $\alpha'(e^B)<0$. The critical points of $F_\lambda$ satisfy
\[ (\epsilon^2\alpha'(s)+((1-\lambda)+\lambda\alpha'(s)) \eta^2 f(p))ds+((1-\lambda)s+\lambda\alpha(s)) \eta^2 df(p)=0.\]
For $\epsilon>\eta>0$ sufficiently small, it follows that critical points $(t,p)$ of $F_\lambda$ are in bijection with critical points $p$ of $f$, and have $t$-coordinate satisfying $A \le t < B$. Moreover, since $\alpha''(s)<0$ for $s \in [e^A,e^B]$, the critical points of $F_\lambda$ are all non-degenerate.

Recall that $e^A$ is the unique critical point of $\alpha(s)$, and that $\alpha(e^A)>0$. All critical points of $F_1$ are thus contained in the hypersurface $\{t = A\}$. Since $F_1$ is of the form $\alpha(s)(\epsilon^2+\eta^2 f)$ in a neighbourhood of this hypersurface, we may assume that $(F_1,dt \otimes dt+g)$ is a Morse-Smale pair for which
\[(C^\bullet_{\mathrm{Morse}}(F_1),d_{F_1})=(C^{\bullet-1}_{\mathrm{Morse}}(f),d_f)\]
under the obvious identification of critical points.

Consider the canonical projection
\[ \pi \co U=[0,B] \times \Lambda \to \Lambda,\]
which maps gradient flow lines of $(F_\lambda,dt \otimes dt+g)$ to (reparametrised) gradient flow lines of $(f,g)$.

Any gradient flow line of $F_\lambda$ which is tangent to the $t$-direction in $U$ must be contained entirely in a line of the form $\R \times \{p\}$, where $p$ is a critical point of $f$. In particular, such a gradient flow line cannot connect two critical points.

Consequently, a non-trivial gradient flow line of negative expected dimension occurring in the family $(F_\lambda,dt \otimes dt+g)$, which connects two critical points, projects under $\pi$ to a non-trivial (reparametrised) gradient flow-line of $(f,g)$ in $\Lambda$, which moreover connects the corresponding critical points of $f$. An index computation implies that they both have the same expected dimension. 

Since $(f,g)$ is a Morse-Smale pair by assumption, the above argument shows that there cannot be any non-trivial gradient flow lines of negative expected dimension in the family $(F_\lambda,dt \otimes dt+g)$. It follows that this family induces a trivial cobordism of the rigid gradient flow lines.

After a generic perturbation of this family, we conclude that
\[(C^{\bullet}_{\mathrm{Morse}}(F_0),d_{F_0})=(C^{\bullet}_{\mathrm{Morse}}(F_1),d_{F_1}),\]
where the complex on the left is obtained by some generic perturbation of the pair $(F_0,dt \otimes dt+g)$.

Finally, this shows that we may assume that the sub-complex
\[(C^{\bullet-1}_{\mathrm{Morse}}(f),d_f) \subset (C^\bullet_{\mathrm{Morse}}(F_-),d_{F_-})\]
has a differential given by
\[d_{F_-}|_{C(f)}=d_{F_0}|_{C(f)}=d_{F_1}|_{C(f)}=d_f.\]
\end{proof}

\bibliography{references}
\bibliographystyle{QT}

\end{document}